\documentclass[12pt,twoside]{amsart}

\usepackage{lipsum,amssymb,latexsym,bbm,graphicx,epsfig,epic,eepic,oldgerm,psfrag,comment}
\usepackage{amsmath,amsfonts,amscd,mathrsfs,amsthm,nicefrac}
\usepackage[dvipsnames]{xcolor}
\definecolor{amber}{rgb}{1.0, 0.75, 0.0}

\usepackage{tikz,wasysym}     
\usepackage{marvosym}

\theoremstyle{plain}

\newtheorem{theorem}{Theorem}
\newtheorem{lemma}[theorem]{Lemma}
\newtheorem{proposition}[theorem]{Proposition}

\newtheorem*{remark*}{Remark}
\newtheorem*{remarks*}{Remarks}

\newtheorem*{example*}{Example}
\newtheorem*{examples*}{Examples}

\newtheorem*{definition*}{Definition}

\newtheorem*{question*}{Question}

\newtheorem*{claim*}{Claim}

\newtheorem{openproblem}{Open Problem}

\def\1{\:\!}
\def\2{\;\!}
\def\s{\smallskip}
\def\m{\medskip}

\def\ga{\alpha}

\def\gg{\gamma}
\def\gd{\delta}

\def\gve{\varepsilon}

\def\go{\omega}

\def\gs{\sigma}

\def\cm{{\mathcal M}}
\def\ct{{\mathcal T}}

\def\Imm{\operatorname{Imm}} 
\def\Mono{\operatorname{Mono}}
\def\Int{\operatorname{Int}}  
\def\GL{\operatorname{GL}}

\def\Z{\operatorname{Z}}
\def\reg{\operatorname{reg}}
\def\Ham{\operatorname{Ham}}
\def\Lag{\operatorname{Lag}}
\def\st{\operatorname{st}}

\def\CC{\mathbb{C}}
\def\DD{\mathbb{D}}

\def\RR{\mathbb{R}}

\def\ZZ{\mathbb{Z}}

\def\C{\mathbb{C}}

\def\R{\mathbb{R}}

\def\Z{\mathbb{Z}}

\def\RP{\operatorname{\mathbb{R}P}}
\def\CP{\operatorname{\mathbb{C}P}}

\def\inn{\operatorname{in}}
\def\out{\operatorname{out}}
\def\area{\operatorname{area}}
\def\SL{\operatorname{SL}}

\def\pp{\partial}

\def\ni{\noindent}
\def\b{\bigskip}
\def\m{\medskip}

\def\id{\mbox{id}}

\def\proof{\noindent {\it Proof. \;}}

\newcommand{\proofend}{\hspace*{\fill} $\square$\\}
\newcommand{\diam}{\hspace*{\fill} $\Diamond$\\}
\newcommand*\circled[1]{\tikz[baseline=(char.base)]{
            \node[shape=circle,draw,inner sep=2pt] (char) {#1};}}

\usepackage{amstext, amsthm, amscd, mathrsfs}
\usepackage{enumitem}
\usetikzlibrary{matrix}

\newcommand{\OP}{\operatorname}

\newcommand{\pr}{\operatorname{pr}}

\newtheoremstyle{TheoremNum}
    {\topsep}{\topsep}              
    {\itshape}                      
    {}                              
    {\bfseries}                     
    {.}                             
    { }                             
    {\thmname{#1}\thmnote{ \bfseries #3}}
\theoremstyle{TheoremNum}

\theoremstyle{definition}

\theoremstyle{remark}

\makeatletter \let\@wraptoccontribs\wraptoccontribs \makeatother

\begin{document}

\title{Lagrangian knots and unknots -- an essay}

\author{Leonid Polterovich} 
\author{Felix Schlenk}


\thanks{GDR is supported by the Knut and Alice Wallenberg Foundation under 
the grants KAW~2021.0191 and KAW~2021.0300, 
and by the Swedish Research Council under the grant number~2020-04426.}

\address{(G.\ Dimitroglou Rizell) 
Department of Mathematics,
Uppsala University, Box 480, SE-751 06 Uppsala, Sweden}
\email{georgios.dimitroglou@math.uu.se}

\address{(L.\ Polterovich)
School of Mathematical Sciences,
Tel Aviv University,
Tel Aviv 69978, Israel}
\email{polterov@post.tau.ac.il}
\thanks{LP is partially supported by the Israel Science Foundation grant 1102/20}

\address{(F.\ Schlenk)
Institut de Math\'ematiques,
Universit\'e de Neuch\^atel,
Rue \'Emile Argand~11,
2000 Neuch\^atel,
Switzerland}
\email{schlenk@unine.ch}


\date{\today}

\begin{abstract}
In this essay dedicated to Yakov Eliashberg we survey the current state of the field of Lagrangian
(un)knots, reviewing some constructions and obstructions along with a number of unsolved questions.
The appendix by Dimitroglou Rizell provides a new take on local Lagrangian knots.
\end{abstract}

\maketitle

\vspace{-0.7cm}
\begin{center}
{\small with an appendix by Georgios Dimitroglou Rizell}
\end{center}

\tableofcontents

Yasha Eliashberg has made several inventions in symplectic topology
that opened a new door and that in his and other researcher's hands developed to rich theories.
In this paper we illustrate this in just one example.
We weave a narrative around Lagrangian knots with several goals in mind: 
to honor Yasha's contributions, 
to survey the current state of the field along with some unresolved questions, 
and to describe a selection of constructions and arguments within the field 
that lend themselves to brief and straightforward exposition, 
apt for a celebratory essay.

\section{Background, basics, and motivation}

\subsection*{Background}

The celebrated Eliashberg--Gromov 
$C^0$-rigidity theorem from \cite{El87} and~\cite[Section 3.4.4]{Gr86}
states that the group of symplectic diffeomorphisms
of a symplectic manifold is $C^0$-closed in the group of all diffeomorphisms.
This rigidity result established the existence of symplectic topology.
On the other hand, there were known flexibility phenomena for symplectic
manifolds due to Gromov's $h$-principles. The co-existence of rigidity
and flexibility in symplectic topology was of great interest ever since,
see Eliashberg's informative survey~\cite{El15}.

One subject in symplectic geometry in which rigid and flexible features go hand in hand
is the problem of Lagrangian knots and unknots,
formulated by Arnold~\cite[$\S$ 6]{Ar86} in~1986.
We have chosen this topic since we find it beautiful and since there has been much recent progress
on Lagrangian knots.

\subsection*{Basic notions}

A submanifold $L$ of a symplectic manifold~$(M,\omega)$ is called {\it Lagrangian}\/ if
$\omega$ vanishes on~$TL$ and if $L$ has half the dimension of~$M$.
Examples are the zero section of cotangent bundles~$T^*L$, endowed with their
canonical symplectic form $\sum_j dp_j \wedge dq_j$.
Furthermore, the product of Lagrangian submanifolds is Lagrangian in the product
of their ambient symplectic manifolds.
In the standard symplectic vector space $(\R^{2n}, \omega_0)$ we therefore have
Lagrangian tori (products of circles in the symplectic coordinate planes),
and through Darboux charts Lagrangian tori exist in all symplectic manifolds.

Roughly speaking, a Lagrangian is knotted if it cannot be deformed to a standard
model Lagrangian.
In $T^*L$ the model is the zero section, and in~$\RR^{2n}$ the models are the product tori.
In the complex projective plane~$\CP^2$, that we always endow with a Fubini--Study form,
the model is the Clifford torus or the real projective plane,
and in a product of two 2-spheres the model is the product of the equators or
(if the spheres have equal area) the anti-diagonal.

Of course, one must specify what ``deform" means. We distinguish several equivalence relations.
Assume that $L,L'$ are Lagrangian sub\-mani\-folds in a symplectic manifold~$(M,\omega)$.
We say that $L,L'$ are

\m \ni
{\bf homologous}\,
if they represent the same homology class in $H_n(M;\ZZ)$.

\s \ni
{\bf smoothly isotopic}\,
if there exists a smooth path of submanifolds from $L$ to~$L'$.

\s \ni
{\bf Lagrangian isotopic}\,
if there exists a smooth path of Lagrangian submanifolds from $L$ to~$L'$.

\s \ni
{\bf Hamiltonian isotopic}\,
if there exists a Hamiltonian isotopy $\phi_t$ of $(M,\omega)$
with $\phi_0 = \id$ and $\phi_1(L) = L'$.

\s
In Sections 1 and 2 we will look at Lagrangian unknots and knots with respect to 
these equivalence relations. 
In Section~3 we look at Lagrangian embeddings that are not regularly homotopic through 
Lagrangian immersions.

One could also ask that $L$ and $L'$ come with a parametrization, and that the above
equivalence relations take into account the parametrization,
but we will only look at unparametrized Lagrangians.

\subsection*{Motivation}

Lagrangian submanifolds arise in symplectic topology on several occasions:
invariant tori of classical mechanics, real parts of complex manifolds equipped with a K\"{a}hler form, and graphs of symplectic diffeomorphisms.
The problem of Lagrangian knots was to some extent motivated by these areas.

To illustrate the connection with mathematical physics,
recall that in~\cite[$\S$~9]{Ar86} Arnold emphasized
(in a dadaist manner typical for him) that
``such natural problems and theorems of symplectic topology as the problem of Lagrangian 
knots... were discovered only as a result of experiments in laser optics and the analysis of the variational principles of Percival, Aubri, and others, connected with the theory of corrosion."

The counterpart of Lagrangian submanifolds in complex analysis are {\it totally real}\/ submanifolds~$L$
of K\"{a}hler manifolds $(M,J)$, i.e., those submanifolds~$L$ for which $JT_xL$ is transverse to~$T_xL$ for every $x \in L$.
Given a K\"{a}hler form~$\omega$ with associated Riemannian metric $\omega(\xi,J\eta)$ on~$M$,
the Lagrangian submanifolds are characterized by the fact that
$JT_xL$ is {\it orthogonal}\/ to~$T_xL$.
If we require the angle between $JT_xL$ and~$T_xL$ to lie in $(\pi/2-\varepsilon, \pi/2]$,
then the corresponding class of submanifolds satisfies an $h$-principle, i.e., is flexible:
Any two such ``$\varepsilon$-Lagrangian embeddings" are isotopic through $\varepsilon$-Lagrangian embeddings,
see~\cite{Gr86} 
\footnote{In our interactions (LP), Yasha often referred to Gromov's book \cite{Gr86}. 
On one occasion, when I couldn't find a statement Yasha claimed was there, 
I asked him to point out the page. Yasha did so promptly, but the statement wasn't there. 
When I gently asked Yasha for a precise location, he responded, 'Between the lines!'}
and also the Intrigue~E1 in the seminal book~\cite{CEM24} by Cieliebak, Eliashberg and Mishachev.
Thus rigidity can (and does!) appear only for Lagrangian submanifolds, i.e.,
when the angle is~$\pi/2$ everywhere.

Finally, if we think on Lagrangians as generalized morphisms in the symplectic category,
the counterpart of the problem of Lagrangian knots is the description of the
symplectic mapping class group, i.e.\ the group formed by the set of connected
components of the symplectomorphism group. 
It is worth mentioning here that the square of Seidel's famous Dehn twist of~$T^*S^2$ 
takes each fiber~$F$ of the cotangent bundle to an exotic ``fiber"~$F'$, 
that coincides with~$F$ outside a compact subset, 
and is isotopic to~$F$ by a smooth compactly supported isotopy which, however, 
cannot be made Hamiltonian~\cite{Se99}.

While being formulated in the early heroic period of symplectic to\-po\-lo\-gy,
the problem of Lagrangian knots remains widely open until now. 
With these motivations at hand,
let us discuss some results and open problems.

\section{Unknots}

The first result on Lagrangian unknots is Eliashberg's theorem from~\cite{El95}
on Lagrangian cylinders:
In the standard $(\RR^4, \omega_0)$ with coordinates $x_1,y_1,x_2,y_2$
consider the straight Lagrangian cylinder
$$
L_0 \,=\, \{ x_1^2+y_1^2 = 1,\, x_2=0,\, 0 \leq y_2 \leq 1 \}
$$
connecting the two hyperplanes
$\Pi_j = \{ y_2 = j\}$, $j=0,1$.
Let $U$ be the region bounded by $\Pi_0$ and~$\Pi_1$.

\begin{theorem}
Every Lagrangian cylinder $L$ in $U$ that agrees with $L_0$ near the boundary is Hamiltonian isotopic to~$L_0$
inside~$U$.
\end{theorem}

\ni
{\it Rough outline of a proof.}
The key technique in the proof is filling with pseudo-holomorphic discs ($J$-discs, for short),
that had been developed by Gromov~\cite{Gr85} and Eliashberg~\cite{El90}.
A filling of the cylinder~$L$ by $J$-discs is a smooth family 
$D_{t \in [0,1]}$
of compact, embedded, and
disjoint $J$-holomorphic discs whose boundary circles foliate~$L$.
%
%
For the straight cylinder~$L_0$, the existence of a filling for the standard~$J$ is obvious.
If one can find a filling of~$L$, it is clear that $L_0$ and~$L$ are smoothly isotopic,
and a few more beautiful and elementary geometric arguments in~\cite[$\S$ 2.5]{El95},
some of which we repeat below, 
imply that one can find a Lagrangian isotopy in~$\RR^4$.
Since the symplectic areas of all non-contractible closed embedded curves on the two cylinders are~$\pi$,
this isotopy can then be included into a Hamiltonian isotopy.
These arguments also show that if the filling discs~$D_t$ lie in~$U$, 
then the Hamiltonian isotopy can be taken with support in~$U$.
While this property of the~$D_t$ is not clear for the filling constructed in~\cite{El95},
it does hold for a filling obtained by neck stretching, as explained to us
by Eliashberg:
Choose $A$ so large that the non-standard Lagrangian cylinder~$L$ is contained in 
$\{ (z_1,z_2) \mid \pi |z_1|^2 \leq A \}$, where $z_j = x_j+iy_j$.
Then compactify the disc of area~$A$ in $\{ (z_1,0) \}$ to~$\CP^1$.
So $\CP^1$ is divided into two discs, a small one of area~$1$, and a big one of area~$A-1$.
Now we do neck stretching along
the boundary of a little disc bundle~$D^*_{\varepsilon} L$  in~$T^*L$,
as in \cite[$\S$\21.3--1.4]{EGH00} or \cite[$\S$\22]{CM18}.
Endow the cylinder~$L$ with the flat metric, so that we know the closed geodesics.
With the right choice of the neck stretching family~$J_s$, each holomorphic sphere in the class~$[\CP^1]$ 
in $\CP^1 \times \CC$
breaks into an inner disc in~$D^*_{\varepsilon} L$ asymptotic to a closed geodesic,
a straight cylinder in the symplectization part, and an outer disc.
We also choose~$J_s$ such that the unit discs in $\{ (z_1,0,j) \} \subset \Pi_j$, $j=0,1$, 
are holomorphic.
Then compactifications of the inner discs form a filling of~$L$ by $J$-discs contained in~$U$.

\begin{figure}[h]
 \begin{center}
  \psfrag{L}{$L$}
  \psfrag{L0}{$L_0$}
  \psfrag{W1}{$W_1$}
  \psfrag{We}{\textcolor{red}{$W_\gve$}}
  \psfrag{f}{$\widetilde \varphi$}
  \leavevmode\epsfbox{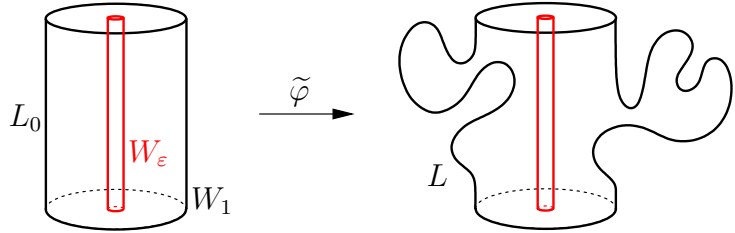}
 \end{center}
 \caption{The symplectomorphism $\widetilde \varphi$ from the filling of $L_0$ to the filling of~$L$.}
 \label{fig-LL}
\end{figure}

Denote the given Lagrangian embedding 
$$
L_0 = \left\{ |z_1| = 1, \, x_2=0,\, 0 \leq y_2 \leq 1 \right\} \,\to\, L
$$
by $\varphi$. 
From this filling one constructs a symplectic embedding $\widetilde \varphi \colon W_1 \to U$
that extends $\varphi$ to the full cylinder
$$ 
W_1 = \left\{ |z_1| \leq 1, \, x_2=0,\, 0 \leq y_2 \leq 1 \right\} \;\:
$$
filling $L_0$
and that is the identity near the discs $\{ y_2 =0\}$ and~$\{ y_2=1 \}$.
After applying a Hamiltonian isotopy supported in~$U$, 
we can assume that $\widetilde \varphi$ is the identity on a thin full cylinder
$$ 
W_{\gve} = \left\{ |z_1| \leq \gve, \, x_2=0,\, 0 \leq y_2 \leq 1 \right\} .
$$
Now take for $t \in [\gve, 1]$ the rescaling
$$
c_t \colon (x_1,y_1,x_2,y_2) \mapsto (t \1 x_1, t \1 y_1, t^2 \1 x_2, y_2) .
$$
It is a diffeomorphism of $U = \{ 0 \leq y_2 \leq 1\}$ such that 
$c_t^*  \omega_0 = t^2 \1 \omega_0$. 
The maps 
$$
\widetilde \varphi_t \,:=\, c_t^{-1} \circ \widetilde \varphi \circ c_t, \quad t \in [\varepsilon, 1]
$$
are symplectic embeddings $W_1 \to U$
that are the identity near the discs $\{ y_2 =0\}$ and~$\{ y_2=1 \}$ and such that 
$\widetilde \varphi_\gve = \id$ and $\widetilde \varphi_1 = \widetilde \varphi$.
Restricting to $\partial W_1 = L_0$ we obtain the sought-after Lagrangian isotopy from $L_0$ to~$L$
that fixes a neighbourhood of the boundary and is supported in~$U$.

\b
For every Lagrangian cylinder~$L$ as above define its ``unknotting size"~$x_2 (L)$ 
as the smallest number~$s$ such that there exists a Hamiltonian isotopy from~$L$ to~$L_0$ 
with support in $U \cap \{|x_2| \leq s \}$.
Note that $x_2(L) \geq \min \{ s \mid L \subset \{ |x_2| \leq s\}\}$
and that $x_2(L) =0$ only if $L=L_0$.

\begin{openproblem} 
What can be said about $x_2(L)$\,?
\end{openproblem}

In \cite{El95}, Eliashberg constructed the filling of $L$ by first deforming the symplectic form~$\omega_0$
to one for which $L$ becomes symplectic.
This idea was used many times later on and in particular in
the proof of the following theorem of Eliashberg and Polterovich~\cite{EP93} 
and also in the proof of Theorem~\ref{t:nolocal} below.

\begin{theorem}[Smooth Unknottedness in $T^*S^2$ and $T^*T^2$]\
\label{t:S2T2}

\s \ni
Assume that $\Sigma$ is the 2-sphere or the 2-torus.
Then every embedded Lagrangian surface $L \subset T^*\Sigma$
that is homologous to the zero section
is smoothly isotopic to the zero section.
\end{theorem}

This theorem was improved much later:
Note that every Lagrangian isotopy of a 2-sphere extends to a Hamiltonian isotopy,
while this is not so for tori.
A Lagrangian $L \subset T^*L$ is called exact
if the restriction of the Liouville form $\sum_j p_j \2 dq_j$ to~$L$ is exact.
Such Lagrangians have a chance to be Hamiltonian isotopic to the zero-section,
and in fact Arnold had already conjectured in~\cite{Ar86}
that this is always the case.
Staying with surfaces here, we refer to~\cite{AbKr18} for the so far best general result 
on this 
`Nearby Lagrangian conjecture' in any dimension. 

\begin{theorem}[Hamiltonian uniqueness in $T^*S^2$ and $T^*T^2$]\ \label{t:Hamunique}

\s \ni
{\rm (i)}
Every Lagrangian sphere in $T^*S^2$ is Hamiltonian isotopic to the zero-section
(Hind~\cite{Hi12}).

\s \ni
{\rm (ii)}
Every exact Lagrangian torus in $T^*T^2$ is Hamiltonian isotopic to the zero-section
(Dimitroglou Rizell, Goodman, and Ivrii~\cite{DGI16}).
\end{theorem}

Hamiltonian uniqueness also holds for $\RP^2$ in~$T^*\RP^2$
(Hind, Pinsonnault and Wu~\cite{HPW16},
as well as Adaloglou~\cite{Ada22}), 
and similarly every Lagrangian sphere in the product of two 2-spheres of equal area is
Hamiltonian isotopic to the anti-diagonal (Hind~\cite{Hi04}).
The proofs of these results all use the stretching the neck technique
and the compactness theorem in symplectic field theory from
Eliashberg, Givental, and Hofer~\cite{EGH00},
and Bourgeois, Eliashberg, Hofer, Wysocki, and Zehnder~\cite{BEHWZ03}.

\m
The first attempt to create a Lagrangian knot is certainly to take a Lagrangian~$L$
and to modify it near a point. This cannot work in view of the following result of
Eliashberg and Polterovich~\cite{EP96}.

\begin{theorem}[Local unknottedness] 
\label{t:nolocal}
Let $L$ be a Lagran\-gian submanifold of $\RR^4$ that is diffeomorphic to
a plane and coincides with a Lagrangian plane~$L_0$ outside a compact subset.
Then there exists a compactly supported Hamiltonian isotopy taking $L$ to~$L_0$.
\end{theorem}

We here give an outline of the proof from~\cite{EP96}. 
Another argument can be found in the appendix by Dimitroglou Rizell.

\begin{figure}[htb]
 \begin{center}
  \psfrag{D0}{$\Delta_0$}
  \psfrag{1}{$\circled{I}$}
	\psfrag{2}{$\circled{II}$}
	\psfrag{3}{$\circled{III}$}
	\psfrag{4}{$\circled{IV}$}
  \psfrag{5}{$\circled{V}$}
	\psfrag{C}{\mbox{cylinder $C$}}
	\psfrag{K}{\textcolor{red}{\mbox{knot $L \supset \Delta$}}}
  \psfrag{Q}{\mbox{hypersurface $Q$}}
  \psfrag{L}{$L$}
  \psfrag{L0}{\textcolor{red}{\mbox{$L_0 \supset \Delta_0$}}}
  \psfrag{zo}{\mbox{zoom out}}
	\psfrag{ppp}{\textcolor{ForestGreen}{\mbox{plug $F_0$}}}
  \psfrag{F}{\textcolor{ForestGreen}{\mbox{$F$}}}
	\psfrag{WWW}{\mbox{Weinstein map}}
  \psfrag{le}{\textcolor{blue}{\mbox{lens}}}
	\psfrag{E}{$E_{\inn}$}
  \psfrag{l0}{$\ell_0$}
  \psfrag{lt}{$\ell_t$}
	\psfrag{l1}{$\ell_1$}
	\psfrag{L=}{{\footnotesize $L = \Gamma (\ell_0)$}}
	\psfrag{Lt}{{\footnotesize $L_t = \Gamma (\ell_t)$}}
	\psfrag{L1}{{\footnotesize $L_1 = \Gamma (\ell_1)$}}
  \psfrag{hy}{\textcolor{blue}{\mbox{hypersurface } $Q$}}
  \psfrag{fil}{\mbox{disc filling}}
	\psfrag{fi}{\mbox{final isotopy}}
	\psfrag{cbg}{\mbox{correcting by gluing}}
  \leavevmode\epsfbox{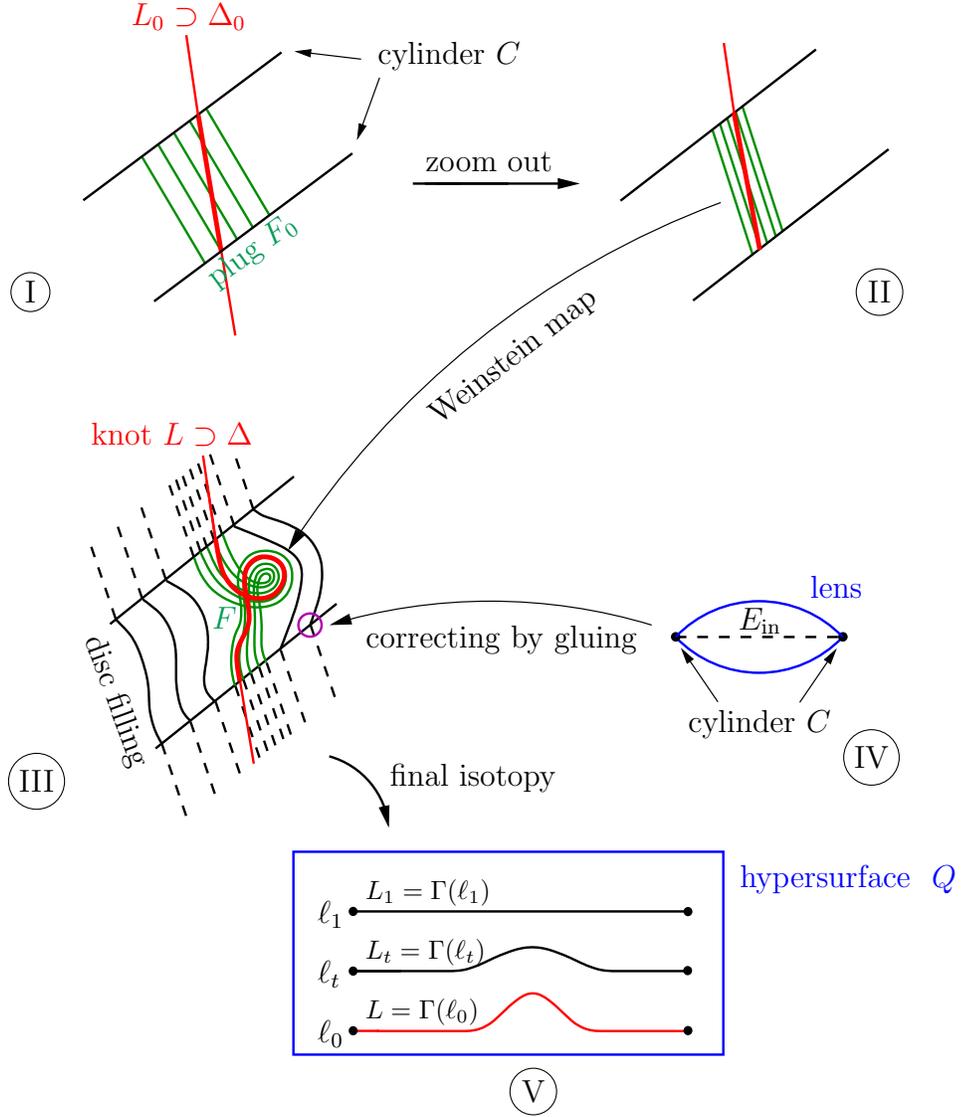}
 \end{center}
 \caption{Local unknottedness - comics}
 \label{comic}
\end{figure}

The main difficulty is that one cannot fill by discs a Lagrangian plane.
To resolve it, we introduce an auxiliary object, a totally real cylinder~$C$ lying
in a hyperplane $E \subset \R^4=\C^2$. 
This cylinder is part of the following local model. 
The cylinder~$C$ splits~$E$ into an inner solid cylinder~$E_{\inn}$ and its complement~$E_{\out}$. 
Fix a Lagrangian plane $L_0 \subset E$ that intersects $E_{\inn}$ in a disc~$\Delta_0$.
There exists a ``plug" $F_0 = D^2 \times (-\varepsilon,\varepsilon) \subset E_{\inn}$ which contains
$\Delta_0$ and such that the discs $D^2 \times \{\text{point}\}$ are affine and symplectic.
The plug is plotted on Figure~\ref{comic}.I, where the discs appear in green. 
For a zoomed out version see Figure~\ref{comic}.II.

Assume now that $L$ is a Lagrangian knot coinciding with~$L_0$ outside a neighbourhood of the
center of~$\Delta_0$. The cylinder~$C$ splits $L$ into a knotted sub-disc~$\Delta$,
and its complement, which coincides with $L_0 \cap E_{\out}$. 
By Weinstein's Lagrangian tubular neighbourhood theorem we can ``insert" the plug, 
now called~$F$, by a symplectic embedding to a tiny neighbourhood of~$\Delta$ (see Figure~\ref{comic}.III).  
By the filling by pseudo-holomorphic discs technique
(applied to a carefully chosen almost complex structure on~$\R^4$)
we extend the foliation on~$F$ by symplectic discs to a filling of the whole cylinder~$C$.

Next we design a pseudo-convex lens (see Figure~\ref{comic}.IV) whose faces are suitable deformations of~$E_{\inn}$ 
with fixed boundary. 
Its interior contains the discs from the filling, and  
the pseudo-convexity provides a control on the geometry of their boundaries. 
This enables us to glue these discs with a natural family of punctured symplectic planes in~$E_{\out}$,
yielding a family of symplectic manifolds diffeomorphic to~$\R^2$.

These symplectic planes foliate a hypersurface $Q \subset \R^4$ containing the knot~$L$. 
Furthermore, this foliation is standard at infinity. The characteristics of~$Q$ are transverse 
to the planes of the foliation, and hence have a very simple dynamics. 
We unknot $L$ by using this dynamics as follows. 
For a curve $\ell \subset Q$ denote by $\Gamma(\ell)$ the union of characteristics passing through~$\ell$ 
(see Figure~\ref{comic}.V).
One easily represents the knot~$L$ as $\Gamma(\ell_0)$ for some line~$\ell_0 \subset Q$. 
Now move $\ell_0$ in the direction ``orthogonal" to~$L$ to a line~$\ell_1$ which is sufficiently 
remote from~$\ell_0$. This family of lines $\ell_t$ gives a desired Lagrangian isotopy $\Gamma (\ell_t)$ 
of~$L$ to an (affine!) Lagrangian 
plane $L_1 = \Gamma(\ell_1)$.  
This completes the outline of the proof of Theorem~\ref{t:nolocal}.

\medskip
Cot\'e and Dimitroglou-Rizell \cite{CDR23} extended Theorem~\ref{t:nolocal} 
to cotangent fibers over open surfaces of finite type. 
Recall, however, that the analogous statement fails for cotangent fibers over~$S^2$
by Seidel's theorem from~\cite{Se99}.

\begin{openproblem} \label{op:surfaces}
Does the analogue of Theorem~\ref{t:nolocal} hold for cotangent fibers over closed surfaces of positive genus?
\end{openproblem}

\medskip
The following problem is wide open.

\begin{openproblem} \label{op:higher}
Does Theorem~\ref{t:nolocal} extend to higher dimensions?
\end{openproblem}

For the weaker result that $L$ as in Theorem~\ref{t:nolocal} is {\it smoothly}\/ isotopic to~$L_0$
a much easier proof can be given, see~\cite{EP97}.
In higher dimensions $n \geq 4$, a smooth isotopy from $L$ to~$L_0$ exists
since then any two embeddings of $S^n$ into~$S^{2n}$ are smoothly isotopic, \cite{Ha61}.
But there do exist embeddings of $S^3$ into~$S^6$ that are not isotopic to the standard embedding, 
\cite{Ha62},
and we do not know if one of these embeddings can be used
to give a negative answer to Open Problem~\ref{op:higher}.

\medskip
A somewhat less local version of Lagrangian knot in~$\RR^4$ would be a compact knot in~$\RR^4$, that one could
then map into other symplectic manifolds by Darboux charts.
The torus is the only closed orientable  surface that admits Lagrangian embeddings into~$\RR^4$.
The first assertion of the following result, again from~\cite{DGI16}, shows that
up to Lagrangian isotopy such local knots do not exist either.

\begin{theorem} \label{t:Lagiso}
{\rm (i)}
All Lagrangian tori in $\RR^4$ are Lagrangian isotopic.

\s
{\rm (ii)}
The same holds true in $\CP^2$ and in the product $S^2 \times S^2$
of spheres of equal area.
\end{theorem}

\section{Knots}

By Theorem \ref{t:nolocal}, Lagrangian knots can only be obtained by some global construction.
By now, Lagrangian knots are known ``at all levels":
homologous but not smoothly isotopic,
smoothly isotopic but not Lagrangian isotopic,
Lagrangian isotopic but not Hamiltonian isotopic:

The following results show that Theorems~\ref{t:S2T2} and~\ref{t:Hamunique} do not hold 
for all symplectic 4-manifolds.

\begin{theorem}
{\rm (i)}
There exist simply connected symplectic 4-manifolds
(for instance the plumbing of two copies of~$T^*S^2$)
that contain infinitely many smoothly isotopic Lagrangian spheres which are pairwise not Lagrangian isotopic
(Seidel~\cite{Se99}).

\s
{\rm (ii)}
There exist simply connected symplectic 4-manifolds
that contain infinitely many homologous Lagrangian tori which are pairwise not smoothly isotopic
(Vidussi~\cite{Vi06}).
\end{theorem}

By Theorem~\ref{t:Lagiso} (i), the only possibility left for a compact orientable 
Lagrangian knot in~$\RR^4$
is a torus that is not Hamiltonian isotopic to a product torus.
Such a torus was constructed by Chekanov in~\cite{Ch96}.
His example was the first Lagrangian knot found.
We here describe the Chekanov torus in the form given by Eliashberg and Polterovich in~\cite{EP97}, 
see also Chekanov and~Schlenk~\cite{CS10}.
Take any closed curve $\gamma (t)$ in $\RR^2$ that is contained in an open half plane 
and encloses area~$1$.
Then view~$\RR^2$ as the complex diagonal of~$\CC^2$, and sweep the curve~$\gamma$
by the anti-diagonal $S^1$-action,
\begin{equation} \label{def:Ch}
\Theta \,:=\, \left\{ \frac{1}{\sqrt 2} \, \bigl( e^{2\pi i \vartheta} \, \gamma (t), e^{- 2\pi i \vartheta} \, \gamma (t) \bigr) \right\} .
\end{equation}
Then $\Theta$ is a Lagrangian torus in $\RR^4$.
We will see later that $\Theta$ does not depend on the precise choice of the curve~$\gamma$.

Let us recall the notion of monotonicity:
For every Lagrangian submanifold $L \subset (M,\omega)$
there are two homomorphisms on the relative homotopy group~$\pi_2(M,L)$,
namely the area class, that measures the symplectic area of a representing disc,
and the Maslov class, that measures how much $TL$ twists along the boundary of a disc
relative to the symplectic structure.
For instance, the Maslov index (i.e., the value of the Maslov class)
of a disc in~$\CC$ is~$2$, since the tangent lines to the boundary circle
make two full turns relative to a fixed direction if we go around the boundary circle once.
A Lagrangian submanifold is called {\it monotone}\/ if these two homomorphisms are positively proportional.
Intuitively, monotone Lagrangians are ``symmetric":
there is a balance between the symplectic size and the symplectic twisting.
For instance, a product torus $T(a,b) \subset \CC^2$ consisting of two circles enclosing area $a$ and~$b$
is monotone iff $a=b$. And $\Theta$ is also monotone.

In \cite{Ch96},
the torus $\Theta$ was distinguished from the monotone product torus $\mathbf{T} := T(1,1)$ of equal area class by looking at the values of a suitable symplectic capacity at nearby tori.
In \cite{EP97}, the distinction was done by the count of $J$-holomorphic discs:
Given a closed Lagrangian $L \subset (M,\omega)$, 
fix a point~$p$ on~$L$, choose an $\omega$-compatible almost complex structure~$J$
on~$M$, and let $\nu (L)$ be the number of $J$-holomorphic discs
of Maslov index~2 with boundary on~$L$ and passing through~$p$.
If $L$ is monotone, then Gromov's compactness theorem implies that $\nu (L)$ does not depend on
the generic choice of $p$ and~$J$.
Furthermore, $\nu$ is invariant under Hamiltonian isotopies.

We have $\nu (\Theta) = 1$ but $\nu (\mathbf{T}) = 2$.
Indeed, for the standard complex structure on $\CC^2$,
the disc in the complex diagonal with boundary $\{ (\gamma(t), \gamma (t))\}$ is unique,
and for the product torus there are only the two coordinate discs.
This count of Maslov index~2 holomorphic discs developed into a powerful invariant for Lagrangian submanifolds,
the so-called disc potential.

To understand the difference between $\mathbf{T}$ and~$\Theta$ more geometrically, 
we come back to the construction of~$\Theta$ in~\cite{EP97},
where $\Theta$ was included into a non-standard fibration of~$\RR^4$ whose fibers are almost
all Lagrangian tori.
See also Section~7.1 on the Auroux system in Evans' book~\cite{Ev23}.
Take the Hamiltonian 
\begin{equation} \label{e:H}
H \colon \CC^2 \to \RR, \quad (z_1,z_2) \mapsto \pi \bigl( |z_1|^2-|z_2|^2 \bigr) 
\end{equation}
as well as the map 
\begin{equation} \label{e:F}
F \colon \CC^2 \to \CC, \quad (z_1,z_2) \mapsto z_1 z_2 .
\end{equation}
For every $a \in \RR$ and every simple closed curve $\gamma \subset \CC$ set
\begin{equation} \label{e:Tag}
T_a(\gamma) \,:=\, \left\{ (z_1,z_2) \mid H(z_1,z_2) = a, \, F(z_1,z_2) \in \gamma \right\} .
\end{equation}
These sets are all invariant under the $S^1$-action 
$$
(z_1,z_2) \mapsto (e^{2\pi i \vartheta}z_1, e^{-2\pi i \vartheta} z_2)
$$
generated by~$H$.
If $a \neq 0$ or if $0 \notin \gamma$, then $T_a(\gamma)$ is a Lagrangian torus,
but if $a=0$ and $0 \in \gamma$, then $T_a(\gamma)$ is an immersed Lagrangian $2$-sphere (the Whitney sphere), 
that can be visualized as a ``croissant" or as a pinched torus.
If we foliate~$\CC$ by nested loops and one point, then we obtain a foliation of 
$\CC^2$ by Lagrangian tori, one pinched torus, and a cylinder over the point. 
This ``almost toric fibration" was taken up later by Auroux in~\cite{Au07} in the context of Mirror symmetry, 
and we shall encounter it soon again for the construction of exotic tori in~$\CP^2$.

Take $a=0$ and look at curves $\gamma$ enclosing area~$\frac{1}{\pi}$.  
If $\gamma_1, \gamma_2$ are such curves going around the origin~$0$,
we can take a Hamiltonian isotopy moving $\gamma_1$ to~$\gamma_2$ in~$\CC \setminus \{0\}$,
and then use the $S^1$-action to lift this isotopy to a Hamiltonian isotopy of~$\RR^4$
that moves $T_0(\gamma_1)$ to $T_0(\gamma_2)$.
Since for $\gamma$ the circle centred at~$0$ we have $T_0(\gamma) = \mathbf{T}$,
all such tori are Hamiltonian isotopic to the Clifford torus~$\mathbf{T}$.
Similarly, any such curve~$\gamma$ not going around~$0$ yields a torus~$T_0(\gamma)$ 
that is Hamiltonian isotopic to the Chekanov torus~$\Theta$.
On the other hand, if $\gamma_1$ goes around the origin while $\gamma_2$ does not, 
then a Hamiltonian isotopy from $\gamma_1$ to~$\gamma_2$ must cross~$0$ 
and hence does not lift to a Hamiltonian isotopy from $T_0(\gamma_1)$ to~$T_0(\gamma_2)$:
The lifted tori degenerate to a pinched torus along the deformation and hence we do not obtain an isotopy of embedded Lagrangian tori.
This is in accordance with the fact proven above that there is no Hamiltonian isotopy 
from $T_0(\gamma_1)$ to~$T_0(\gamma_2)$ at all.

It is still an open question whether in $\RR^4$ the Chekanov torus is the only monotone Lagrangian knot:

\begin{openproblem} \label{p:R4}
Is every monotone Lagrangian torus in $\RR^4$ Hamiltonian isotopic to (a scaling of) the Clifford torus
or the Chekanov torus?
\end{openproblem}

\begin{figure}[h]
 \begin{center}
  \psfrag{16}{$1$}
  \psfrag{13}{$2$}
  \psfrag{12}{$3$}
	\psfrag{m1}{$\mu_1$}
	\psfrag{m2}{$\mu_2$}
  \leavevmode\epsfbox{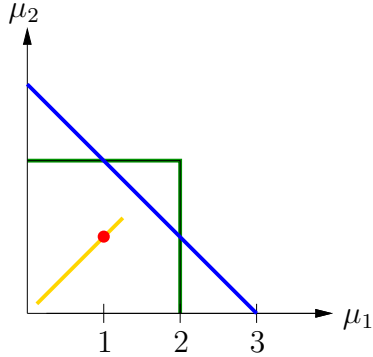}
 \end{center}
 \caption{The images under $\mu$ of $\Theta$ and $\mathbf{T}$ in 
$B^4(3) \cap D(2) \times D(2)$.}
 \label{fig.moment}
\end{figure}

The Chekanov torus $\Theta$ naturally sits as a monotone torus in $\CP^2$ and in~$S^2 \times S^2$:
To have it monotone we must scale the Fubini--Study form to integrate to~$3$ over the complex line, 
and we must take spheres of area~$2$.
These manifolds are compactifications of the open ball $B^4(3)$ of capacity $\pi r^2=3$
and of the product of two discs~$D(2)$ of area~$2$.
If we choose the curve $\gamma$ enclosing area~$1$ 
in the construction of~$\Theta$ in the half-disc of area~$\frac 32$,
then indeed $\Theta \subset B^4(3) \cap D(2) \times D(2)$.
This is best seen on the image of the moment map 
$$
\mu \colon \CC^2 \to \RR_{\geq 0}^2, \quad 
      \mu (z_1,z_2) = \left( \pi |z_1|^2, \pi |z_2|^2 \right)
$$
that generates the standard Hamiltonian torus action
\begin{equation} \label{e:T2}
(z_1,z_2) \,\mapsto\, \left( e^{2\pi i\vartheta_1} z_1, e^{2\pi i \vartheta_2} z_2 \right) ,
\end{equation}
see Figure~\ref{fig.moment}.
These monotone Lagrangian tori in $\CP^2$ and $S^2 \times S^2$ are, again, not 
Hamiltonian isotopic to the model tori (i.e., the Clifford torus and the product of equators), 
see~\cite{CS10}.
In contrast to~$\RR^4$, however, 
in these compact manifolds many more Hamiltonian isotopy
classes of monotone Lagrangian tori are known! 
The following result should also be compared with Theorem~\ref{t:Lagiso}~(ii).

\begin{theorem} \label{t:vianna}
There are infinitely many Hamiltonian isotopy classes of monotone Lagrangian tori in~$\CP^2$
and in the product $S^2 \times S^2$ of spheres of equal area.
\end{theorem}

These tori were predicted by Galkin and Usnich~\cite{GU10}
by a mirror-symmetry argument, and they  
were independently constructed by Vianna \cite{Vi16,Vi17} and by Galkin--Mikhalkin~\cite{GM22}
in two different ways:
As the central fibers of almost toric fibrations, and
as the tori obtained from the central fiber of certain weighted projective planes after a smoothing.
We briefly outline both approaches. 

\m \ni
{\bf Almost toric fibrations.}
The projective plane
$\CP^2$ carries a ``best Lagrangian fibration":
Viewing $\CP^2$ as the compactification of the open ball~$B^4(3)$ by a sphere~$\CP^1$ at infinity,
the Hamiltonian torus action~\eqref{e:T2} extends to~$\CP^2$.
The image of $\CP^2$ under the extended moment map is the closed triangle~$\Delta$ 
shown in Figure~\ref{fig-poly}.I below. 
The fiber over an interior point is a Lagrangian 2-torus, while the fibers over the edges and vertices
are circles and points, respectively.
The only monotone fiber is the one over the red centre point, it is the Clifford torus.

To find knotted monotone Lagrangian tori, we look at ``worse" fibrations,
which have more complicated singularities than the toric fibration.
An almost toric fibration (ATF for short) of~$\CP^2$ is the ``next best" kind of fibration:
Some of its fibres are immersed Lagrangian spheres with one double point,
i.e.\ pinched tori like the above $T_0 (\gamma)$ for a simple closed curve $\gamma$ through the origin.
$\CP^2$ admits many~ATFs. 
The simplest one has one pinched torus but one point fibre less than the standard fibration.
The passage to this new fibration is called `nodal trade' and 
was first described by Symington in~\cite{Sy03} and then in great detail in Evans' monograph~\cite{Ev23}. 
We here give a more direct description.

\begin{figure}[h]
 \begin{center}
  \psfrag{r1}{$\rho_1$}
  \psfrag{r2}{$\rho_2$}
  \psfrag{e}{$\varepsilon$}
  \psfrag{d}{$\delta$}
	\psfrag{Ld}{$L_\delta$}
  \psfrag{rr}{$\rho_1\rho_2 = \pi^2 r^2$}
  \psfrag{F}{$F=z_1z_2$}
  \psfrag{0}{$0$}
  \psfrag{p}{$p$}
	\psfrag{g0}{\textcolor{blue}{$\gamma_0$}}
	\psfrag{DC}{$D(r)$}
  \psfrag{U}{$\mu(U)$}
	\psfrag{L}{$L_{\gve}$}
  \leavevmode\epsfbox{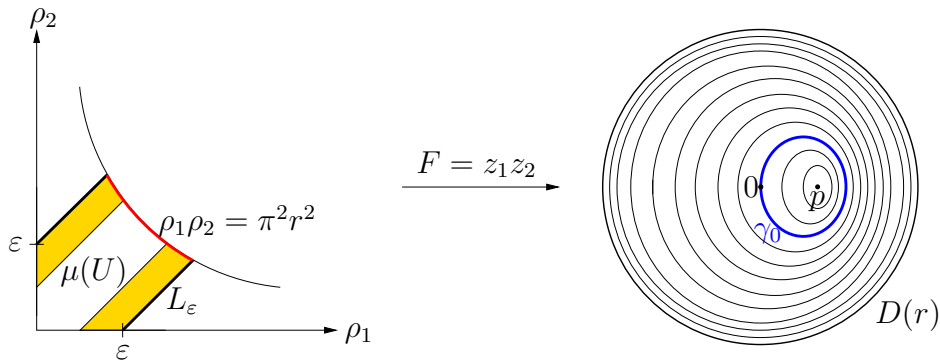}
 \end{center}
 \caption{The region $\mu(U)$ and the foliation of $D(r)$.}
 \label{fig-foliation}
\end{figure}

\begin{figure}[htb]
 \begin{center}
  \psfrag{v1}{$v_1$}
  \psfrag{v2}{$v_2$}
	\psfrag{v3}{$v_3$}
  \psfrag{v1'}{$v_1'$}
  \psfrag{v2'}{$v_2'$}
	\psfrag{v3'}{$v_3'$}
  \psfrag{I}{I}
  \psfrag{II}{II}
  \psfrag{III}{III}
  \psfrag{IV}{IV}
  \psfrag{D}{$\Delta$}
	\psfrag{D'}{$\Delta'$}
	\psfrag{Dp}{$\Delta_+$}
  \psfrag{Dm}{$\Delta_-$}
	\psfrag{R}{$R$}
  \psfrag{S}{$S$}
	\psfrag{R'}{$R'$}
	\psfrag{SD}{$S \Delta_-$}
  \leavevmode\epsfbox{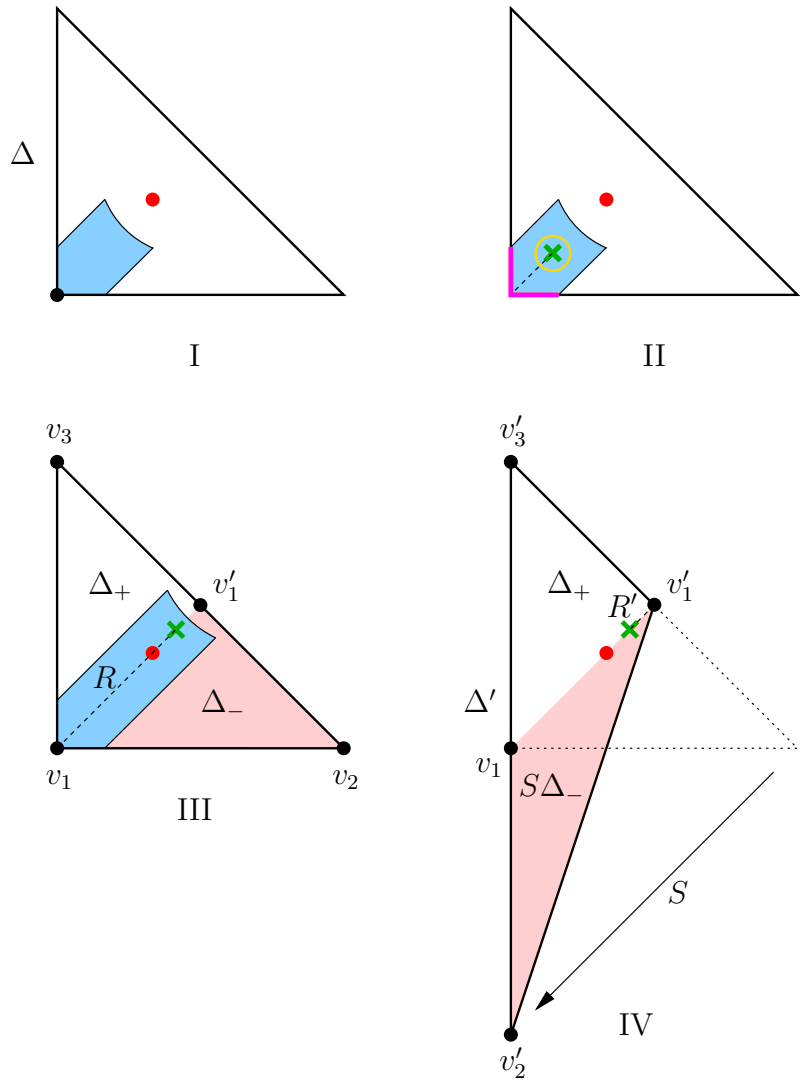}
 \end{center}
 \caption{From the toric fibration of $\CP^2$ (I) to an almost toric fibration (II), 
and the first mutation (III and~IV).}
 \label{fig-poly}
\end{figure}

We take up the ATF of~$\RR^4$ constructed above, 
and restrict it to the 4-ball 
\begin{eqnarray*}
U &=&
\left\{ 
(z_1,z_2) \in \CC^2 \mid - \gve \leq H(z_1,z_2) \leq \gve, \, |F(z_1,z_2)| \leq r
\right\} \\
&=& 
\left\{ 
(z_1,z_2) \in \CC^2 \mid - \gve \leq \rho_1-\rho_2 \leq \gve, \, \rho_1 \rho_2 \leq  \pi^2 r^2
\right\}
\end{eqnarray*}
where $H$ and $F$ are given by \eqref{e:H} and \eqref{e:F},
where $\gve, r >0$ and where we use action variables $\rho_j = \pi |z_j|^2$.
We now foliate the disc $D(r) \subset \CC$ of radius~$r$ by simple closed curves and the point~$p \neq 0$
as in Figure~\ref{fig-foliation}.
Near the boundary, the leaves are concentric circles.
This yields an ATF on~$U$ whose fibers are 
the Lagrangian tori~$T_a(\gamma)$ defined in~\eqref{e:Tag} with $-\gve \leq a \leq \gve$, 
one pinched torus~$T_0(\gamma_0)$,
and the (interval of) circles over~$p$.
As an aside, we note that by a result of Dimitroglou Rizell~\cite{Rizell2019}, 
any closed embedded or immersed Lagrangian in~$U \setminus F^{-1}(p)$ with the same classical invariants as a fiber 
is Hamiltonian isotopic to a fiber. 

Since~$H = \rho_1-\rho_2$ is constant on each fiber $T_a(\gamma)$, each fibre lies over 
a segment $L_\delta := \{\rho_1-\rho_2 = \delta \} \cap \mu(U)$ parallel to the diagonal. 
Since the leaves of the foliation near the boundary of~$D(r)$ are concentric circles,
over a neighbourhood of the red curve $\{ \rho_1 \rho_2 = \pi^2 r^2\} \cap \mu(U)$
our ATF agrees with the standard toric fibration of~$\R^4$ by product tori.
In order to glue the entire ATF with the standard toric fibration outside~$U$, 
we modify the ATF over the two yellow bands.
For $\delta \in [\frac 12 \gve, \gve]$ consider the set
$\mu^{-1}(L_\delta)$, that is the union of tori $T_\delta (\gamma)$ and of one circle (over~$p$).
The $S^1$-action restricts to a free $S^1$-action on~$\mu^{-1}(L_\delta)$,
and the reduced space $\mu^{-1}(L_\delta) / S^1$ is a closed disc~$D_\delta$
that is foliated in two ways:
by concentric circles and one point (coming from the standard toric fibration restricted to~$L_\delta$)
and by the loops $T_\delta (\gamma) / S^1$ and one point.
This second foliation~${\mathcal F}_\delta$ also consists of nested embedded loops and one point, 
and the two foliations agree near the boundary of~$D_\delta$.

Let $\varphi_\delta$ be a smooth family of Hamiltonian diffeomorphisms of~$D_\delta$ with compact support
in the interior, such that $\varphi_\delta = \id$ for $\delta \in [\frac 12 \gve, \frac 23 \gve]$ and such that
$\varphi_\delta$ takes the foliation~${\mathcal F}_\delta$ to the concentric foliation
for $\delta \in [\frac 34 \gve, \gve]$.
Lifting the image foliation of~$\varphi_\delta$ by the $S^1$-action we obtain a new Lagrangian 
foliation of~$\mu^{-1}(L_\delta)$.
Applying the same construction over the yellow band near~$L_{-\gve}$
we obtain an ATF on~$U$ that smoothly fits with the standard toric fibration on
the complement of~$U$.

\m
This change of the fibration, the `nodal trade', is somewhat schematically represented by Figure~\ref{fig-poly}.II: 
Outside the blue region~$\mu(U)$, the fibration II agrees with the standard fibration
in~I.
Over the blue region, the lower left point fiber in~I is traded against the node
(double point) of the pinched torus, that lies over the green cross. 
The fibers over the two purple segments are circles. 
These are the circles that (before the above interpolation)
were the circles over the point~$p$ in Figure~\ref{fig-foliation}, 
and contrary to what II suggests, the union of these circle fibers is a {\it smooth}\/ cylinder. 
All other fibers over the blue region are Lagrangian 2-tori. 
The dashed line indicates that the topological monodromy of the $T^2$-bundle over a circle
around the cross is non-trivial. 

Recall that the region $U = U_{\gve,r}$ in the above construction of the nodal trade
depends on parameters $\gve,r>0$.
We now take $r$ so large that the base point of the pinched torus (the cross) lies 
on the other side of the central point.
The new ATF is shown in~Figure~5.III.
Let's call this seemingly harmless operation a `nodal slide'.
For the central torus this operation is less harmless: 
The torus over the red point in~II is the Clifford torus of~$\CP^2$,
while by our construction the torus over the red point in~III
is the torus in~$\CP^2$ obtained as above from the Chekanov torus
$\Theta \subset B^4(3) \subset \CP^2$,
and we already know that this torus is not Hamiltonian isotopic to the Clifford torus.

To bring into play Markov numbers later on, we now describe the ATF~III 
of~$\CP^2$ in a different way.
Cut the triangle~$\Delta$ along the diagonal line spanned by the dashed ray~$R$ into 
the upper and lower triangles $\Delta_+$ and~$\Delta_-$.
Let $S$ be the linear transformation that is the identity on~the diagonal line and shears down~$\Delta_-$
as shown in Figure~\ref{fig-poly}.IV:
$S$ fixes $v_1$ and~$v_1'$ and takes $v_2$ to~$v_2'$.
This yields a new triangle $\Delta' := \Delta_+ \cup S \Delta_-$.
Coming back to~III, let~$\Delta_{\reg}$ be the regular part of~$\Delta$, 
namely the interior of~$\Delta$ deprived from the point at the cross.  
Recall that the monodromy of the $T^2$-bundle over~$\Delta_{\reg}$ around the cross is not trivial. 
In fact, we can construct this $T^2$-bundle by taking the trivial~$T^2$-bundle
over~$\Delta_{\reg} \setminus R$ and gluing the $T^2$-fibers along~$R$ by the transpose of the inverse of~$S$.
The same ATF is described by~IV, where the fibrations over~$\Delta_+$ and~$S\Delta_-$ are glued as follows: 
The non-trivial gluing of $T^2$-fibers along~$R$ in~III is now undone
and instead the non-trivial gluing is done over 
the short dashed ray~$R'$ on the other side of the cross.
This passage from III to~IV is called `transferring the cut'. 
While the fibrations in~II and~III are different, those in III and~IV are isomorphic, 
they are just described over different base triangles.

We apply the nodal trade operation described in I $\to$ II also at the upper and lower vertex of~$\Delta'$.
We can now do the same geometric mutation 
(nodal slide followed by transferring the cut) 
at any of the three vertices and then look at the torus over the red point. 
If we do a mutation at the newly created vertex~$v_1'$ in~IV,
we just undo the previous mutation, and obtain again the Clifford torus. 
But if we do a mutation at~$v_2'$ or $v_3'=v_3$, we potentially obtain new tori. 
At any rate, geometric mutations produce a trivalent tree~$\ct$ of triangles with tori
over the red central point:
The vertices of~$\ct$ are the triangles thus obtained, and two vertices of~$\ct$ are connected by an edge
if and only if the two corresponding triangles are related by one geometric mutation.

This tree of triangles $\ct$ can be matched with the Markov tree~$\cm$, whose vertices are labeled by 
triples~$(a,b,c)$ of natural numbers that solve the Markov equation
\begin{equation}
\label{e:markov}
a^2 + b^2 + c^2 \,=\, 3abc.
\end{equation}
Given any such solution, three other solutions can be obtained by keeping two numbers and
replacing the third one, say $a$, by $3bc-a$.
Starting from $(1,1,1)$, these algebraic mutations create the trivalent tree of solutions whose 
beginning is shown below,
and every solution of~\eqref{e:markov} appears in this tree.

\begin{center}
\begin{tikzpicture}[level distance=1cm,
  level 3/.style={sibling distance=6cm},
  level 4/.style={sibling distance=3cm}]
  \node {$(1,1,1)$}
  child { node {$(2,1,1)$}
  child {node {$(5,2,1)$}
    child {node {$(13,5,1)$}
      child {node {$(34,13,1)$}}
      child {node {$(194,13,5)$}}
    }
    child {node {$(29,5,2)$}
    child {node {$(433,29,5)$}}
      child {node {$(169,29,2)$}}
    }}};
\end{tikzpicture}
\end{center}

By now, we have found a matching $\ct \to \cm$
between the tree of triangles~$\ct$ and the Markov tree~$\cm$. 
Denote by $\Delta (a,b,c)$ the triangle in~$\ct$ corresponding to the 
Markov triple~$(a,b,c) \in \cm$.
We can upgrade the matching by relating the geometry of~$\Delta (a,b,c)$ with the Markov triple~$(a,b,c)$.
Given a planar triangle with edges of rational slope, 
we define the weight of a vertex~$v$ as
the absolute value of the determinant of the two primitive lattice vectors parallel to the edges meeting at~$v$. Weights are integers.
For instance, the weights of~$\Delta = \Delta (1,1,1)$ are~$(1,1,1)$ and 
those of~$\Delta' = \Delta (2,1,1)$ are~$(4,1,1)$.
More generally, we have:

\begin{proposition} \label{p:abc}
The weights of the triangle~$\Delta (a,b,c)$ are $(a^2,b^2,c^2)$.
\end{proposition}

\proof
This was shown in~\cite{Vi16}, see also \cite[Appendix~H.2]{Ev23}.
We give a more elementary proof.

\begin{figure}[h]
 \begin{center}
  \psfrag{q1}{$q_1$}
  \psfrag{q2}{$q_2$}
	\psfrag{q3}{$q_3$}
  \psfrag{v1}{$v_1$}
  \psfrag{v2}{$v_2$}
	\psfrag{v3}{$v_3$}
	\psfrag{vh1}{$\hat v_1$}
  \psfrag{vh2}{$\hat v_2$}
	\psfrag{vh3}{$\hat v_3$}
	\psfrag{v1'}{$v_1'$}
  \psfrag{v2'}{$v_2'$}
  \leavevmode\epsfbox{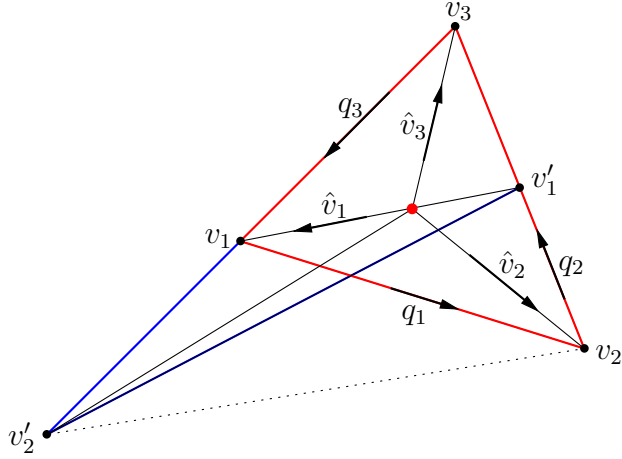}
 \end{center}
 \caption{Primitive vectors associated with $\Delta (a,b,c)$, and the new vertices $v_1'$ and $v_2'$ of $\Delta (a',b,c)$.}
 \label{fig-vq}
\end{figure}

After a translation we can assume that the red central point of $\Delta (a,b,c)$
is at the origin of~$\RR^2$.
For $i=1,2,3$ we denote by~$\hat v_i$ the primitive vector in the direction of the vertex~$v_i$
and by $q_i$ the primitive vector in the direction of the $i$'th oriented edge, 
see Figure~\ref{fig-vq}.
We use cyclic notation mod~$3$ for the indices, e.g.\ $q_{i-1} = q_3$ for $i=1$.
The weight at~$v_i$ is the determinant $w_i := |q_{i-1} \times q_i|$.

Write $(n_1,n_2,n_3) = (a,b,c)$, and consider the assertions

\s
\begin{itemize}
\item[1)]
$w_i = n_i^2$ \, for $i=1,2,3$,

\s
\item[2)]
$\hat v_i = \displaystyle \frac{q_{i-1}-q_i}{n_i}$ \, for $i=1,2,3$.
\end{itemize}

\s \ni
Note that 1) is the claim in Proposition~\ref{p:abc}.
We shall prove 1) and~2) for all triangles in~$\ct$ by induction down the tree~$\ct \cong \cm$:
1) and~2) hold at the root $\Delta (1,1,1)$.
We assume that 1) and 2) hold at $\Delta (a,b,c)$, and shall prove them for the two triangles in~$\ct$
right below~$\Delta (a,b,c)$. 
Let $\Delta (a',b',c')$ be one of these triangles. After renaming, we can assume that
$\Delta (a',b',c')$ is obtained from $\Delta (a,b,c)$ by mutation at~$v_1$.
Hence $a' = 3bc-a$ and $b'=b$, $c'=c$.
Recall that $v_1'$ is the point of intersection of the line through~$v_1$ and the origin with the edge~$[v_2,v_3]$, 
that $v_2'$ is the point of intersection of the line through $v_3$ and~$v_1$ with the line through~$v_2$ parallel to~$[v_1,v_1']$, 
and that $\Delta (a',b,c)$ is obtained from $\Delta (a,b,c)$ by applying to the triangle with 
vertices~$v_1, v_2, v_1'$ the shear~$S$ that fixes the line generated by~$v_1$ and takes $v_2$ to the point~$v_2'$ on the line
trough $v_1$ and~$v_3$.

\begin{lemma}
$S$ is the shear by $\hat v_1$, i.e., 
$$
S (p) = \gs_{\hat v_1} (p) := p + (\hat v_1 \times p) \hat v_1 
\quad \mbox{for all } p \in \RR^2 .
$$
\end{lemma}

\proof
Using both inductive assumptions 1) and 2) we compute
\begin{eqnarray*}
\gs_{\hat v_1} (q_1) &=& q_1 + \frac{1}{n_1^2} \bigl( (q_3-q_1) \times q_1 \bigr) (q_3-q_1) \\
&=& q_1 + \frac{w_1}{n_1^2} (q_3-q_1) \,=\, q_3 .
\end{eqnarray*}
Hence $S$ and $\gs_{\hat v_1}$ both fix the line spanned by~$v_1$ pointwise and take the oriented line
trough $v_1, v_2$ to the oriented line through $v_3,v_1$.
Hence $S = \gs_{\hat v_1}$.
\proofend

Since $S$ fixes the primitive vector $\hat v_1$ and takes the primitive vector~$q_1$
to the primitive vector~$q_3$, it takes an oriented basis of~$\ZZ^2$ to an oriented basis of~$\ZZ^2$, 
i.e., $S \in \SL(2;\ZZ)$.
We can now prove assertions 1) and~2) for $\Delta (a',b,c)$.

\s
Proof of 1): We have $w_3' = w_3 = n_3^2$, and since $S \in \SL (2;\ZZ)$ we have $w_2' = w_2 = n_2^2$.
The weight $w_1'$ at the newly created vertex~$v_1'$ is
\begin{eqnarray*}
w_1' 
&=& \left| S (-q_2) \times q_2 \right| \\
&=& \left| \gs_{\hat v_1} (q_2) \times q_2 \right| \\
&=& \left| \bigl( q_2 + (\hat v_1 \times q_2) \hat v_1 \bigr) \times q_2 \right| \\
&=& (q_2 \times \hat v_1)^2 .
\end{eqnarray*}
For the integer $n_1' := q_2 \times \hat v_1$ we thus have 
$$
n_1' \,=\, q_2 \times \frac{q_3-q_1}{n_1} \,=\, \frac{w_3+w_2}{n_1} \,\stackrel{(\ast)}{=}\,
3n_2n_3-n_1 \,=\, 3bc-a = a' 
$$
where in $(\ast)$ we used 1) and the Markov equation \eqref{e:markov} for $(a,b,c) = (n_1,n_2,n_3)$.
Hence $w_1' = (n_1')^2  = (a')^2$.

\medskip
Proof of 2):
This is clear for the unchanged vector $\hat v_3$.
For $\hat v_2'$ we must show $\hat v_2' = \frac{q_3-Sq_2}{n_2'}$.
This holds since $w_2=w_2'$ implies $n_2=n_2'$ and since
$$
\hat v_2' \,=\, S \hat v_2 \,\stackrel{2)}{=}\, S \left( \frac{q_1-q_2}{n_2}\right)
\,=\, \frac{1}{n_2} (q_3-Sq_2) . 
$$
It remains to show that $\hat v_1' = \frac{Sq_2-q_2}{n_1'}$. This holds since $\hat v_1' = - \hat v_1$
and 
$$
S q_2 - q_2 \,=\, \gs_{\hat v_1} (q_2) - q_2 \,=\, (\hat v_1 \times q_2) \hat v_1  \,=\, -n_1' \hat v_1.
$$
The induction step is done and Proposition~\ref{p:abc} follows.
\proofend

Once we know that $S \in \SL(2;\ZZ)$, we can 
prove Proposition~\ref{p:abc} in a more direct way:
The integral length of a segment~$s$ in~$\RR^2$ with rational slope can be defined as follows:
Let $A$ be a translation followed by a matrix in $\SL (2,\ZZ)$ that takes $s$ to the $x$-axis. 
Then the integral length of~$s$ is the Euclidean length of $A (s)$. 
Equivalently, the integral length of~$s$ is $|\ell|$ if the difference of the end-points of~$s$ is equal 
to~$\ell  q$ for a primitive vector~$q$.

Let $u_i = \ell_i q_i$ be the oriented edges of $\Delta (a,b,c)$,
so $\ell_i$ is the integral length of~$u_i$.
%
%
After scaling $\Delta = \Delta (1,1,1)$ we can assume that
the edges of $\Delta$ have integral length~$3$.
Since $\Delta (a,b,c)$ is obtained by geometric mutations from~$\Delta$
and since the matrix~$S$ of each geometric mutation belongs to~$\SL (2;\ZZ)$,
the area and integral perimeter of these two triangles are the same:
$$
|u_{i-1} \times u_{i} | = 2 \area \Delta (a,b,c) = 9, \qquad
\ell_1 + \ell_2 + \ell_3 = 9.
$$
The first identity and 
$|u_{i-1} \times u_{i} | = \ell_{i-1} \ell_{i} |q_{i-1} \times q_{i} |
 = \ell_{i-1} \ell_{i} w_i$ yield 
$$
\ell_{i-1} \ell_{i} w_i = 9 .
$$
From these three identities and $\ell_1 + \ell_2 + \ell_3 = 9$
we obtain
\begin{equation} \label{e:M2}
\sqrt{9 w_1 w_2 w_3} \,=\, w_1 + w_2 + w_3.
\end{equation}
Let $\overline \Delta$ be the triangle in~$\ct$ right above~$\Delta (a,b,c)$,
with weights $\overline w_1, \overline w_2, \overline w_3$.
Since $\Delta (a,b,c)$ is obtained from~$\overline \Delta$ by one geometric mutation, two of the weights are equal, say  
$\overline w_2 = w_2$, $\overline w_3 = w_3$.
By induction we assume that the weights of~$\overline \Delta$ are squares. 
Then $w_2 = n_2^2$, $w_3 = n_3^2$, and~\eqref{e:M2} becomes
$$
\sqrt{w_1} \,=\, \frac{w_1 + w_2 + w_3}{3n_2n_3} .
$$
Since the weights $w_i$ are integers, this is a rational number.
Hence $w_1$ is also a square, $w_1 = n_1^2$. 
Now~\eqref{e:M2} reads 
$$
3 n_1 n_2 n_3 \,=\, n_1^2+n_2^2+n_3^2 ,
$$
i.e., $(n_1,n_2,n_3)$ is a Markov triple.
Since a geometric mutation of a triangle in~$\ct$
changes only one weight and an algebraic mutation of a Markov triple
changes only one Markov number, we conclude by induction along~$\ct$
that $(\sqrt{w_1}, \sqrt{w_2}, \sqrt{w_3}) = (n_1,n_2,n_3) = (a,b,c)$.
\diam

Let $T(a,b,c)$ be the torus over the red central point of $\Delta (a,b,c)$. 
To distinguish these tori,
let $\nabla (a,b,c)$ be the convex hull of the set of those elements of~$H_1(T(a,b,c))$ 
that are represented by the boundary of a $J$-holomorphic disc of Maslov index~2 with boundary on~$T(a,b,c)$.
Similar to the counting invariant $\nu (L)$ described early
in this section, the monotonicity of~$T(a,b,c)$ implies that
$\nabla (a,b,c)$ does not depend on the generic choice of~$J$
and is
invariant under Hamiltonian isotopies.
For concreteness we fix a $\ZZ$-basis of $H_1(T(a,b,c))$. 
Then $\nabla (a,b,c)$ becomes a subset of~$\RR^2$.
Let $\Delta^{\!\circ}(a,b,c)$ be the dual of~$\Delta (a,b,c)$, namely the triangle
whose vertices are the primitive lattice vectors outward normal to the edges of~$\Delta (a,b,c)$.
Vianna proved by neck-stretching that 
\begin{equation} \label{e:Via}
\nabla (a,b,c) = A (\Delta^{\!\circ} (a,b,c) ) \quad \mbox{for a matrix $A \in \GL(2;\ZZ)$.}
\end{equation}
Assume now that $T(a,b,c)$ and $T(a',b',c')$ are Hamiltonian isotopic.
Then $\nabla (a,b,c) = \nabla (a',b',c')$.
Hence $\Delta^{\!\circ} (a,b,c) = B (\Delta^{\!\circ} (a',b',c'))$ for a matrix $B \in \GL(2;\ZZ)$, by~\eqref{e:Via}.
Hence $B^T (\Delta (a,b,c)) = \Delta (a',b',c')$.
Hence the weights of $\Delta (a,b,c)$ and $\Delta (a',b',c')$, that are $\GL(2;\ZZ)$-invariants
and equal $\{a^2,b^2,c^2\}$ and $\{a'^2, b'^2, c'^2\}$,
are the same.
Hence $\{a,b,c\} = \{a',b',c'\}$.
This proves Theorem~\ref{t:vianna} for~$\CP^2$. 

\m
Recall from the general Theorem 5~(ii)
that all the tori $T(a,b,c)$ are Lagrangian isotopic. 
This
can be seen directly from their construction.

\begin{proposition} \label{p:Viannaisotpic}
The tori $T(a,b,c)$ are all Lagrangian isotopic.
\end{proposition}

\begin{figure}[h]
 \begin{center}
  \psfrag{q1}{$q_1$}
  \psfrag{q2}{$q_2$}
	\psfrag{q3}{$q_3$}
  \leavevmode\epsfbox{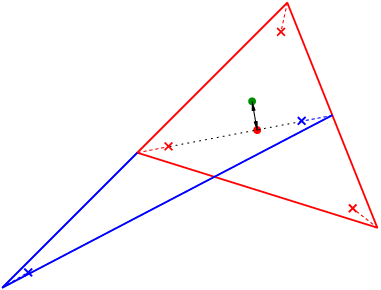}
 \end{center}
 \caption{Proof of Proposition~\ref{fig-proof}.}
 \label{fig-proof}
\end{figure}

\proof
The proof is illustrated in Figure~\ref{fig-proof}.
Consider the torus $T(a,b,c)$ over its red base point in $\Delta (a,b,c)$,
and one of the two tori, say $T(a',b,c)$, just below it in the Markov tree,
that is obtained by geometric mutation.
Move the red base point a bit up to the green point in the interior 
of the triangle~$\Delta_+$, and call the (non-monotone) torus over the green point~$T_{\gve}$.
This move of the base point corresponds to a Lagrangian isotopy of the fibers 
from $T(a,b,c)$ to~$T_{\gve}$
since the ATF is trivial near the segment from the red point to the green point.
We now do the geometric mutation leading to~$T(a',b,c)$.
Since the green point does not lie on the cut line, 
we can take the support of the nodal slide disjoint from $T_{\gve}$.
And since $T_{\gve}$ lies over the interior of~$\Delta_+$, it is also untouched by the
transferring the cut operation (the half-shear defined by $S$).
Finally, take as before a small Lagrangian isotopy moving $T_{\gve}$ to~$T(a',b,c)$.
Composing the to Lagrangian isotopies from $T(a,b,c)$ to $T_{\gve}$ and from $T_{\gve}$ to~$T(a',b,c)$
we obtain a Lagrangian isotopy from $T(a,b,c)$ to~$T(a',b,c)$.
Now use induction over the Markov tree.
\proofend

In the next section we will also look at Lagrangian immersions. 
We already prove here that as Lagrangian immersions, the tori $T(a,b,c)$ are equivalent
in the strongest possible way.
A regular homotopy of Lagrangian immersions $f_t \colon L \to L_t$ in a symplectic manifold $(M,\omega)$
is said to be {\it exact}\/ 
if the area homomorphism $[\omega]_t \colon \pi_2(M,L_t) \to \RR$ does not depend on~$t$.

\begin{proposition} \label{p:GLh}
Assume that $L_0,L_1$ are monotone Lagrangian tori in~$\CP^2$ that are Lagrangian isotopic.
Then $L_0,L_1$ are regularly homotopic through exact Lagrangian immersions.
\end{proposition}
%

\proof
This essentially follows from the Gromov--Lees $h$-principle.
We adapt the proof of 24.3.1 in~\cite{CEM24} to our situation.
Let $f_t \colon L \to L_t$ be a Lagrangian isotopy from $L_0$ to~$L_1$.
Extend $f_t \circ f_0^{-1}$ to a smooth family of diffeomorphisms 
$\widetilde f_t \colon \mathcal{N}(L_0) \to \mathcal{N} (L_t)$
between tubular neighbourhoods $\mathcal{N}(L_t) \subset M$.
Since $L_t$ is Lagrangian, $(\widetilde f_t)^* \omega$ vanishes on~$L_0$.
Hence there exist 1-forms~$\widetilde \ga_t$ on $\mathcal{N}(L_0)$ such that
\begin{equation} \label{e:gat}
(\widetilde f_t)^* \omega = d \widetilde \ga_t \quad \mbox{ and } \quad
\mbox{$\ga_t := \widetilde \ga_t |_{L_0}$  
is closed.}
\end{equation}

Since $L_0$ and $L_1$ are monotone, 
\begin{equation} \label{e:oo01}
[\go]_j \colon \pi_2(\CP^2,L_j) \to \RR \quad \mbox{agree for $j=0,1$}. 
\end{equation}
Choose a basis $\{\gg_k \}$ of $\pi_1(L_0) \cong H_1(L_0;\ZZ)$
and classes $\{ [D_k] \}$ in $\pi_2(\CP^2,L_0)$ such that $\partial [D_k] = \gamma_k$.
Then choose closed 1-forms $\beta_0, \beta_1$ on~$L_0$ such that
\begin{equation} \label{e:jk}
[\ga_j+\beta_j] \gg_k \,=\, [\go]_j [D_k] \quad \mbox{for $j=0,1$ and $k=1,2$.}
\end{equation}
Replace $\widetilde \ga_t$ by $\widetilde \ga_t + t \pr^* \beta_1 + (1-t) \pr^* \beta_0$,
where $\pr \colon \mathcal{N}(L_0) \to L_0$ is the projection, 
and call this form again~$\widetilde \ga_t$.
Then \eqref{e:gat} still holds, and now $[\ga_0] = [\ga_1] \in H^1(L_0;\RR)$ 
by~\eqref{e:oo01} and~\eqref{e:jk}.
Choose $h \colon L_0 \to \RR$ with $\ga_1 = \ga_0 + tdh$, and extend~$h$ to $\widetilde h$ 
on~$\mathcal{N}(L_0)$.
Making a final correction, we replace $\widetilde \ga_t$ by $\widetilde \ga_t - \pr^* \ga_0 - t d\widetilde h$,
and call this form again~$\widetilde \ga_t$.
Then~\eqref{e:gat} still holds, and now $\ga_t := \widetilde \ga_t |_{L_0} = 0$ for $t=0,1$. 
In particular, the inclusion
\begin{equation} \label{e:iota}
\iota \colon L_0 \,\to\, \bigl( \mathcal{N}(L_0), d \widetilde \alpha_t \bigr), \quad t \in [0,1]
\end{equation}
defines an isotopy of Lagrangian embeddings that for $t = 0,1$ are exact in
the sense that the 1-forms $\iota^* \widetilde \alpha_t = \ga_t$ on $L_0$ are exact.
(In fact, they vanish.)

We can now apply to the Lagrangian isotopy \eqref{e:iota} 
the $h$-principle~24.3.1 in \cite{CEM24}, and obtain a 
regular homotopy of Lagrangian immersions
$$
g_t \colon L_0 \,\to\, g_t(L_0) \subset 
           \bigl( \mathcal{N}(L_0), d \widetilde \alpha_t \bigr), \quad t \in [0,1]
$$
such that $g_0=g_1= \iota$ and such that all $g_t$ are exact, i.e., 
$g_t^* \widetilde \ga_t$ are exact for all $t \in [0,1]$.
The composition 
$$
\widetilde f_t \circ g_t \circ f_0 \colon L \,\to\,  
     (\widetilde f_t \circ g_t) (L_0) =: L_t'
\,\subset\, (M,\omega)
$$ 
is a regular homotopy of Lagrangian immersions from $f_0$ to~$f_1$
for which the area homomorphism
$[\omega]_t \colon \pi_2(M,L_t') \to \RR$ is constant on~$[0,1]$.
Indeed, if $u \colon D \to \CP^2$ is a disc with boundary~$\gamma$
on~$L_0$, 
then for each $t \in [0,1]$ we obtain a disc~$D_t$ with boundary on~$L_t'$
by adding to~$u$ the cylinder 
$$
z_t \colon S^1 \times [0,t] \to \CP^2, \quad
(\theta, s) \mapsto \bigl( \widetilde f_s \circ g_s \bigr) (\gamma(\theta)) ,
$$
and $D$ and $D_t$ have the same $\omega$-area because
$\int_{S^1 \times [0,t]} z_t^* \omega =0$.
\proofend

We finally sketch the approach to Theorem~\ref{t:vianna} of Galkin and Mikhalkin in~\cite{GM22}.

\m \ni
{\bf Tori from degenerations.}
Let $Y \subset \CP^N$ be a complex projective variety with singular set $\Sigma$,
and assume that $\pi \colon Y \to \DD$ 
is a holomorphic submersion on $Y \setminus \Sigma$
that takes $\Sigma$ to~the center~$0$ of the open unit disc $\DD \subset \CC$. 
For $t \in \DD$ set $X_t = \pi^{-1}(t)$.
One says that $X_{t \neq 0}$ is a smoothing of~$X_0$ and that $X_0$ is a degeneration 
of~$X_{t \neq 0}$.
The restriction of the Fubini--Study form~$\omega$ on~$\CP^N$ turns the $X_{t \neq 0}$,
$X_0 \setminus \Sigma$,
and~$Y \setminus \Sigma$ into symplectic manifolds.
The symplectic planes that are $\omega$-orthogonal to $TX_t \subset TY$ define a connection 
and a symplectic parallel transport on $Y \setminus \Sigma$.
It follows that all $X_{t \neq 0}$ are symplectomorphic. 
Moreover, parallel transport takes a Lagrangian torus in $X_0 \setminus \Sigma$ to a 
Lagrangian torus in every~$X_{t \neq 0}$.

Now take $X_0$ to be a weighted projective plane $\CP^2(a^2,b^2,c^2)$.
This is a singular toric variety (with at most three singular points), and the moment map triangle can
be taken to be $\Delta (a,b,c)$.
Hacking and Prokhorov showed in~\cite{HP10} 
that $\CP^2(a^2,b^2,c^2)$ admits a smoothing as above exactly if $(a,b,c)$
solves the Markov equation.
In this case, we obtain from the fiber over the central point of $\Delta (a,b,c)$
a monotone Lagrangian torus $\widetilde T(a,b,c)$ in $X_t = \CP^2$.
Galkin and Mikhalin proved~\eqref{e:Via} for $\widetilde T(a,b,c)$.
Hence the tori $\widetilde T(a,b,c)$ are not Hamiltonian isotopic for different Markov triples,
implying again Theorem~\ref{t:vianna} for~$\CP^2$.

\m
It is widely believed that the answer to the following question is `yes'.

\begin{openproblem} 
{\rm
Is it true that for each Markov triple $(a,b,c)$
the tori $T(a,b,c)$ and $\widetilde T(a,b,c)$ constructed above
are Hamiltonian isotopic\,?
}
\end{openproblem}

This would follow at once from an affirmative answer to

\begin{openproblem} 
{\rm
Is every monotone Lagrangian torus in $\CP^2$ Hamiltonian isotopic to a torus $T (a,b,c)$\,?
}
\end{openproblem}

In a similar way, one can construct infinitely many monotone Lagrangian tori in $S^2 \times S^2$
that are pairwise not Hamiltonian isotopic. 
They are now parametrized by the nodes of a tree that consists of triangles and quadrilaterals, 
see~Vianna~\cite{Vi17}, and Pascaleff and Tonkonog~\cite{PT20}.

\begin{openproblem} 
{\rm
Is every monotone Lagrangian torus in $S^2 \times S^2$ Hamiltonian isotopic to one of these tori\,?
}
\end{openproblem}

\subsection*{Some applications}
The above exotic tori in $\CP^2$ and $S^2 \times S^2$ have various applications
to symplectic topology. 
We just state four of them. 

\m
(i)
By a result of McDuff~\cite{Mc91}, the space of symplectic embeddings of a closed 4-ball into an open 4-ball
is connected. For other domains, the situation can be very different. 
For instance, there are infinitely many symplectic embeddings of 
the closed polydisc $\overline D(1) \times \overline D(1)$
into~$B^4(3)$ that are not symplectically isotopic.
Brendel, Mikhalkin, and Schlenk~\cite{BMS24} show this by using the tori $T(a,b,c)$ with $a=1$.

\m
(ii)
Yet another equivalence relation on Lagrangian submanifolds $L,L' \subset (M,\omega)$
is defined by asking that there exists a symplectomorphism taking $L$ to~$L'$.
Let $\sigma$ be the symplectomorphism of~$S^2 \times S^2$ that exchanges the factors.
Hind, Mikhalkin, and Schlenk show in~\cite{HMS24} that 
for most ATFs of $S^2 \times S^2$ with quadrilateral base
the monotone torus~$L$ over the center
is not Hamiltonian isotopic to~$\sigma (L)$. 

\m
(iii)
Successively lifting the tori $T(a,b,c)$ in $\CP^2$ to Lagrangian tori in $\CP^3, \CP^4, \dots$,
Chanda, Hirschi, and Wang~\cite{CHW23} constructed
infinitely many Hamiltonian isotopy classes of monotone Lagrangian tori in every~$\CP^n$.

\m
(iv)
Consider the standard contact sphere $(S^{2n+1}, \xi_{\st})$ of dimension at least~$5$.
Lifting the Lagrangian tori in~(iii), 
Blakey, Chanda, Sun, and Woodward~\cite{BCSW24}
obtain infinitely many {\it Legendrian}\/ tori in these spheres 
which are not Legendrian isotopic to each other.

\m
More generally, 
let $L$ and $K$ be monotone Lagrangians in a closed
symplectic manifold $(M,\omega)$ with integral symplectic form, 
and assume that $L$ and~$K$ can be lifted to embedded Legendrians $L'$ and~$K'$
of a prequantization $(P,\alpha)$ of~$(M,\omega)$. 
Consider the following statements:

\begin{itemize}
\item[(1)]
$L$ and $K$ are Hamiltonian isotopic;

\s
\item[(2)] 
$L'$ and $K'$ are Legendrian isotopic.
\end{itemize}
Note that (1) yields (2), while (2) in general does not yield (1).
Thus (2) can be considered as a weak form of unknottedness which is not so much explored yet.
This discussion also extends to more general Bohr--Sommerfeld Lagrangians.

\subsection*{In $\RR^{2n}$}
While for odd $n \geq 5$ all monotone Lagrangian tori in $\RR^{2n}$ 
are smoothly isotopic, 
this is not so for even $n \geq 4$, 
as Dimitroglou Rizell and Evans showed in~\cite{DE12}
by using Haefliger--Hirsch theory 
(cf.\ the results by Borrelli~\cite{Bo01} and Nemirovski~\cite{Ne24} 
for the case of $S^k \times S^1$ discussed at the end of~$\S$\24).

In contrast to $\RR^4$ (see Problem~\ref{p:R4}), 
it is known that for $2n \geq 6$ there are infinitely many Hamiltonian isotopy classes 
of monotone Lagrangian tori in~$\RR^{2n}$ with equal area classes.
Infinitely many such tori were first constructed by Auroux~\cite{Au15} who, again, 
distinguished them by the disc potential, and then by Brendel~\cite{Br23} whose construction is 
inspired by the one of the Chekanov torus in~\eqref{def:Ch}.
An iteration procedure based on~\eqref{def:Ch} was used earlier in~\cite{CS10} to 
construct many (though finitely many) such tori.
The relation between these three sets of monotone Lagrangian tori in $\RR^{2n}$ ($2n \geq 6$)
has not been worked out.

\subsection*{Non-monotone tori}
The above Lagrangian knots are all monotone. 
Monotone Lagrangi\-ans can only exist in symplectic manifolds that are monotone themselves 
(meaning that the area class on $\pi_2(M)$ is positively proportional to the first Chern class).
In contrast, sufficiently small tori in~$\R^{2n}$, monotone or not,
can always be embedded into a given symplectic manifold,
and Brendel showed in~\cite{Br23} that for all ``reasonable" symplectic manifolds of dimension at least~$6$
(including closed ones and cotangent bundles),
every open subset contains infinitely many Lagrangian tori
which are pairwise not Hamiltonian isotopic but
are Lagrangian isotopic and have the same area class.
Thus the occurrence of infinitely many Hamiltonian knotted tori is a purely local phenomenon 
in dimensions~$\geq 6$.
Recall, however, Open Problem~\ref{op:higher}.

In dimension 4, it is harder to find non-monotone tori with the same area class that are not Hamiltonian isotopic. 
Examples in~$\CP^2$ and in the monotone~$S^2 \times S^2$ 
were found by Shelukhin, Tonkonog, and Vianna~\cite[$\S$~7]{STV18} 
and by Fukaya, Oh, Ohta, and Ono~\cite{FOOO12}, 
and in many more symplectic 4-manifolds by Brendel, Hauber, and~Schmitz in the recent work~\cite{BHS24}.
To describe exotic non-monotone tori in~$\CP^2$
we return to Figure~\ref{fig-poly}.
Take a point in the triangle in~I that lies in the interior of the segment from the black point
to the red point. 
The two tori over this point in I and in~III have the same area classes, but are not Hamiltonian isotopic.  
The main open problem on non-monotone tori in dimension~4 is

\begin{openproblem} \label{p:R4non-monotone}
Is every non-monotone Lagrangian torus in $\RR^4$ Hamiltonian isotopic to 
a product torus $S^1(a) \times S^1(b)$\,?
\end{openproblem}

{\footnotesize
When Yasha spent the spring of 2022 at the ITS of ETH Z\"urich,
he proposed to me (FS) ``that we come together for a day".
(Of course, it became two days, with overnight at Ada and Yasha's place.)
In a work with Brendel about Lagrangian pinwheels
I was desperately looking for a certain $J$-sphere, that somehow had to
exist, but I already doubted it did.
Asking Yasha for how to get this sphere, he set off a firework, or
rather a cascade of fireworks:
an idea came up ... mmh,
then another one ... ``something is fishy here" while walking to Coop to buy a bottle of wine,
then a longer line of arguments along neck stretching,
showing it really should exist,
then yet a simpler approach through resolutions ...
At the end I used a theorem of Taubes, but the key point was that I got
totally charged by Yasha and convinced that there {\it is}\/ an argument.
}

\section{Knottedness from the Maslov class}
Another flavor of Lagrangian knots comes from the
existence of Lagrangian embeddings which are not regularly homotopic
through Lagrangian immersions.
We consider here the case of Lagrangian embeddings $f \colon S^k \times S^1 \to \RR^{2k+2}$
which are obtained by a Lagrangian surgery of an immersion $S^{k+1} \to \RR^{2k+2}$ with a
single double point, that is transverse.
We assume throughout that $k \geq 2$.
By Gromov's theorem, one of the generators of the first homology has positive symplectic area.
The minimal Maslov number, denoted by~$\mu(f)$, is defined as the Maslov index of this generator.
The absolute value of this number agrees for Lagrangian embeddings which are regularly homotopic
through Lagrangian immersions.

Let's start with the Whitney immersion
$$
w \colon S^{k+1} \to \CC^{k+1}(p+iq), \quad w(x,y) = (1+iy)x
$$
where $S^{k+1} = \{ (x,y) \in \RR^{k+1} \times \RR \mid |x|^2+y^2 = 1\}$.  
(For $k=1$ this is the 2-dimensional pinched torus encountered earlier.)
Recall that there are two different surgeries~\cite{P-surgery}. 
One of these surgeries always leads to a ``standard" Lagrangian embedding~$f_k$ of~$S^k \times S^1$ with $\mu(f_k)=2$.
When $k \geq 3$ is odd, the second surgery gives rise to a Lagrangian embedding $g_k$ of~$S^k \times S^1$ 
with $\mu(g_k)=k+1$.
When $k \geq 2$ is even, however, the second surgery leads to a Lagrangian embedding
of a non-trivial $S^k$-bundle over~$S^1$, and hence becomes irrelevant to our discussion.

As an application of the $h$-principle for Lagrangians with certain conical singularites over loose Legendrians,
proven by Eliashberg and Murphy in~\cite{EM13},
Ekholm, Eliashberg, Murphy, and Smith \cite{EEMS} constructed for every even~$k$
a Lagrangian immersion of $S^{k+1}$ into~$\RR^{2k+2}$ having a unique transverse double point
that is different from the Whitney immersion.
One of the surgeries performed at this point leads to an 
exotic Lagrangian embedding~$h_k$ of~$S^k \times S^1$
with a number of highly unexpected properties.
First, when $k=2$, the Maslov class of~$h_k$ vanishes. This provides a negative answer to a question
by Audin~\cite[p.\ 622]{Audin}, and in fact debunks a once widely accepted belief that the Maslov
class of a closed Lagrangian submanifold in~$\RR^{2n}$ cannot vanish.
Second, for $k \geq 4$, $h_k$ is negatively monotone with $\mu(h_k) = 2-k$.
The construction of the Lagrangian embeddings $h_k$ is not explicit.

\begin{openproblem}
Find an explicit Lagrangian embedding $S^2 \times S^1 \to \RR^6$ with vanishing Maslov class.
\end{openproblem}

In view of their minimal Maslov numbers, $f_k$ and $g_k$, 
$g_k$ and~$h_k$, 
as well as $f_k$ and~$h_k$ with $k \neq 4$,
are not regularly homotopic through Lagrangian immersions. 
We claim that the embeddings $f_4$ and~$h_4$ are not isotopic through Lagrangian embeddings. 
Indeed, assume on the contrary that such an isotopy, say $u_t \colon S^k \times S^1 \to \R^{2k+2}$ 
with $u_0= f_4$ and $u_1 = h_4$,
does exist. Denote by $\lambda_t$ the corresponding symplectic area class, 
and by $\mu$ the (non-vanishing and $t$-independent) Maslov class in 
$H^1(S^k \times S^1,\Z) = \Z$. Note that $\lambda_t = c(t) \1 \mu$ with $c(0) >0$ and $c(1) < 0$ 
(the latter is a consequence of the negative monotonicity). 
Thus, by continuity, there exists $s \in (0,1)$ with $c(s)=0$.
Thus $\lambda_s = 0$, which contradicts Gromov's theorem on the non-existence of exact Lagrangian embeddings. 
The claim follows. 

\begin{openproblem}
Are the embeddings $f_4$ and~$h_4$ regularly homotopic through Lagrangian immersions?
\end{openproblem}

As pointed out by Georgios Dimitroglou--Rizell, the construction 
of the Lagrangian immersion of~$S^5$ underlying the embedding~$h_4$
actually depends on the choice of the formal Lagrangian homotopy class of the Lagrangian disc-cap from~\cite{EM13},
and so $h_k$ and the solution to Open Problem 9 may depend on this choice.

\begin{proposition}  \label{p:Ju}
For all even $k \geq 2$, the embeddings $f_k$ and~$h_k$
are regularly homotopic through smooth immersions.  
\end{proposition}

This follows from the fact that 
{\it two immersions of a closed mani\-fold~$M$ of odd dimension~$n \geq 3$ into $\RR^{2n}$ 
are regularly homotopic through smooth immersions if and only if 
they have the same number of double points mod 2}.
We learned this statement from a remark on \cite[p.\ 81]{Ju}, 
and Andr\'as Juh\'asz explained to us the following proof.
By the Hirsch--Smale $h$-principle the differential gives a weak homotopy equivalence 
$\Imm(M, \R^{2n}) \to \Mono(TM, T\R^{2n})$ from the space of immersions to the space of 
fiberwise injective bundle homomorphisms $TM \to T \RR^{2n}$ covering a
continuous map $M \to \R^{2n}$.
In particular, the obstructions to a regular homotopy between two immersions lie
in $H^i(M; \pi_i(V_n(\R^{2n})))$, 
where $V_n(\R^{2n})$ is the Stiefel manifold of $n$-frames in~$\R^{2n}$. 
As $\pi_i(V_n(\R^{2n})) = 0$ for $i<n$ and is $\Z_2$ for~$i=n$, 
the only obstruction lies in $H^n(M; \Z_2) = \Z_2$. 
So there are at most two regular homotopy classes. 
The number of double points mod~2 is an invariant under regular homotopy 
(consider double point curves of the trace of the regular homotopy in~$\R^{2n} \times [0,1]$). 
By taking a connected sum with an immersion of~$S^n$ with one double point we see that 
both 0 and~1 can be realized as the number of double points mod~2. 
The claim follows.

A much stronger result was recently obtained by Nemirovski~\cite[Cor.\ 1.2]{Ne24}.

\begin{theorem} \label{t:Ne}
All Lagrangian embeddings of $S^k \times S^1$ in~$\CC^{k+1}$
are smoothly isotopic for $k\neq 3$, 
and for $k=3$ there are two smooth isotopy classes of Lagrangian embeddings.
In particular, the embeddings $f_k$ and~$h_k$ are isotopic through smooth embeddings.
For every $k$, an explicit Lagrangian embedding is given by
$$
S^k \times S^1 \to \RR^{k+1} \oplus \RR^{k+1}, \qquad
(x,t) \mapsto \bigl( x+ \tfrac 12 \sin t \, x ,  \tfrac 12 \cos t \, x \bigr) ,
$$
where $S^k$ is viewed as the unit sphere in $\RR^{k+1}$ and $S^1$ as $\RR/2\pi \ZZ$.
\end{theorem}

The proof is based on the Haefliger--Hirsch embedding theorem and results of Skopenkov~\cite{Sko08} for~$k=2$.
For even $k \geq 2$ there are infinitely many smooth isotopy classes of embeddings of $S^k \times S^1$ in~$\CC^{k+1}$,
and for odd $k \geq 3$ there are two such isotopy classes, see~\cite[p.\ 2]{Ne24}.
For $k \neq 3$, Theorem~\ref{t:Ne} is therefore a rigidity theorem.
Such knotted smooth embeddings have first been found by Hudson~\cite{Hu63},
and two more geometric constructions are given by Skopenkov~\cite[$\S$\23]{Sko}.

Yet another construction yielding all smooth isotopy classes of embeddings for $k \geq 3$
and infinitely many for $k=2$ was shown to us by Nemirovski:
Embed $S^k$ in the standard way in $\R^{k+1}$.
Consider its (trivial) normal bundle~$\nu$ in~$\C^{k+1}$ of rank~$k+2$.
We are going to choose a trivial rank~2 sub-bundle of~$\nu$ and then take small
circles in these normal planes. 
This will define a smooth embedding of $S^k \times S^1$ into~$\CC^{k+1}$.
For the choice of the trivial rank 2 sub-bundle,
we first take any section of $\nu$ consisting of unit vectors, for instance $x \mapsto (x,0_{k+1})$.
This choice is unobstructed, since it corresponds to a map from $S^k$ to~$S^{k+1}$
because $\nu$ is trivial of rank~$k+2$.
Choosing the second normal vector field in~$\nu$ then means
choosing a map from $S^k$ to~$S^k$, which gives an integer~$W$.
The Haefliger--Hirsch embedding theorem and Lemma~\ref{le:Wen} below
show that for even $k \geq 4$ the integer~$W$ is a complete invariant
of isotopy classes of smooth parametrized embeddings, 
and that for odd $k \ge 3$ its residue mod 2 is a complete invariant.
One finally checks that reparametrization of the domain~$S^k \times S^1$
for even $k \geq 4$
identifies the isotopy classes for which the degree of the map $S^k \to S^k$ differs only by the sign, 
and for odd $k \geq 3$
does not change the set of isotopy classes.
For $k=2$, this construction also yields infinitely many isotopy classes of embeddings of $S^2 \times S^1$
into~$\CC^3$, but now there are other ones (see again \cite[p.\ 2]{Ne24}).
The isotopy class of the Lagrangian embeddings is the one obtained by
choosing a map $S^k \to S^k$ of degree~$\pm 1$, 
for instance $x \mapsto (0_{k+1},x)$.


We finally compute the Whitney invariant of the above embeddings. 
Referring to \cite{HH63, Sko} for a general discussion of this invariant, 
we restrict to the case $M_k := S^k \times S^1$ with $k \geq 3$.
Set 
$$
\Z_{M_k} := 
\left\{
\begin{array}{ll}
\ZZ  & \mbox{ if \, $k \geq 4$ is even,} \\ [0.1em]
\ZZ_2 &   \mbox{ if \, $k \geq 3$ is odd,}         
\end{array}
\right.
$$
and let $E(S^k\times S^1)$ be the set of isotopy classes of parametrized smooth embeddings $S^k \times S^1 \to \C^{k+1}$.
Recall that the Haefliger--Hirsch embedding theorem from~\cite{HH63} states that the Whitney invariant 
$$
W \colon E(S^k \times S^1) \to 
H_1 (S^k \times S^1;\Z_{M_k}) = \Z_{M_k} 
$$
is a bijection.
According to \cite[$\S$\25]{Sko}, the Whitney invariant $W(f)$ of a smooth embedding $f \colon S^k \times S^1 \to \C^{k+1}$
can be defined as follows:
Choose a reference embedding $e \colon S^k \times S^1 \to \C^{k+1}$. 
After applying an isotopy, we can assume that $e$ and $f$ agree outside of a closed embedded topological $k+1$-ball 
$B \subset S^k \times S^1$. 
Choose a smooth homotopy $F \colon B \times [0,1] \to \CC^{k+1}$ from $e|_B$ to~$f|_B$ which is constant on the boundary of~$B$, 
and consider the set 
$$
f \cap F := (f |_{M_k \setminus B})^{-1} \bigl( F (B \times [0,1]) \bigr) ,
$$ 
and its closure $\overline{f \cap F}$ in $M_k$.
Then $W(f)$ is defined as 
$$
[\overline{f \cap F}] \,\in\, 
H_1(M_k \setminus \Int B, \pp B; \ZZ_{M_k}) \cong 
H_1(M_k; \ZZ_{M_k}) \cong \ZZ_{M_k} .
$$
Note that $W$ depends on the choice of the reference embedding $e$.
For every $n \in \ZZ$ explicit representatives $e_n$ of the isotopy classes of the
embeddings $S^k \times S^1 \to \C^{k+1}$ described in the above construction of Nemirovski
are given by
\begin{equation}
e_n(x,t) \,=\, \bigl( x + c \2 \cos t \, x, c \2 \sin t \, w_n(x)  \bigr) \,\in\, \R^{k+1} \times \R^{k+1} ,
\end{equation}
where $w_n \colon S^k \to S^k$ is any smooth map of degree~$n$ and where $c \in (0,1)$, so that $e_n$ is an embedding.
We write $W(e_n)$ for the Whitney invariant of~$e_n$ with respect to the reference embedding~$e_0$.

\begin{lemma} \label{le:Wen}
$W(e_n) \,=\,
\left\{ 
\begin{array}{ll}
n \in \ZZ & \mbox{ if $k \geq 4$ is even,} \\  [0.1em]
n \! \mod 2 \in \ZZ_2 & \mbox{ if $k \geq 3$ is odd.} 
\end{array}
\right.
$
\end{lemma}

\proof
Fix $c \in (0,1)$ and $n \in \ZZ$, and consider the two maps $e, \tilde e_n \colon S^k \times S^1 \to \R^{k+1} \times \R^{k+1}$
given by
\begin{eqnarray*}
e (x,t) &=& \bigl( x + c \2 \cos (t+2 \delta) \, x  ,\;  c \, \widetilde \sin \, t \, w_0  \bigr) \\
\tilde e_n (x,t) &=& \bigl( x + c \2 \cos (t+2 \delta) \,x ,\;  c \, \widetilde \sin \, t \, w_n(x) \bigr) , 
\end{eqnarray*}
where $\delta >0$ is small and $\widetilde \sin \, t$ is an odd function on the circle that vanishes exactly on 
$[-\delta, \delta]$ and at~$\pi$, has non-negative derivative on $[-2\delta, 2\delta]$, 
and agrees with $\sin t$ outside $[-2\delta, 2\delta]$.
Furthermore, $w_0$ is a point very close but different from the south pole $p_s$ of~$S^k$
(and in any case such that $x_{k+1}(-w_0) \geq \frac 34$), and $w_n$ is a map of degree~$n$.
Then $e$ and $\tilde e_n$ are embeddings, and $\tilde e_n$ is isotopic to~$e_n$.

To choose $w_n$ specifically, we first take a map $v_n \colon S^k \to S^k$ that takes the southern hemisphere to~$p_s$,
on $\{ x_{k+1} \geq \frac 34\}$ has the form
$$
v_n(x_1,x_2,x_3, \dots, x_{k+1}) \,=\, (r^n e^{in\varphi},x_3, \dots, x_{k+1} ),
$$
where $(x_1,x_2) = r e^{i\varphi}$,
and is such that the last coordinate $x_{k+1}$ is mapped monotone increasingly under~$v_n$.
Then we take $w_n \colon S^k \to S^k$ to be a small deformation of~$v_n$ that agrees with~$v_n$ on $\{ x_{k+1} \geq \frac 34\}$
and maps the southern hemisphere to $w_0$.

The maps $e$ and $\tilde e_n$ agree outside of the topological ball
$$
B = \{ (x,t) \mid x_{k+1} \geq 0, \, \delta \leq t \leq 2\pi-\delta \},
$$
as illustrated in the figure.

\begin{figure}[h]
 \begin{center}
  \psfrag{S}{$S^k$}
  \psfrag{S1}{$t \in S^1$}
	\psfrag{R}{$B$}
	\psfrag{gd}{$\gd$}
	\psfrag{-gd}{$2\pi -\gd$}
  \psfrag{x}{$x_{k+1}=0$}
  \leavevmode\epsfbox{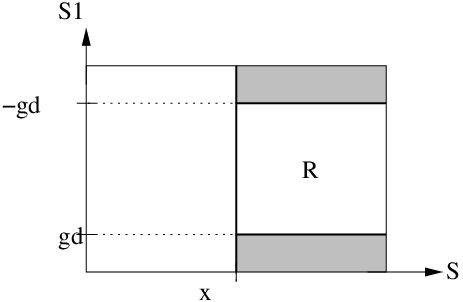}
 \end{center}
 \caption{}
 \label{fig-}
\end{figure}

As homotopy between $e$ and $\tilde e_n$ we take the linear one:
$$
F_n (x,t,s) = 
\Bigl( x + c \2 \cos (t+2 \delta) \, x ,\; 
       c \, \widetilde \sin t \, \bigl( s w_0 + (1-s) w_n(x) \bigr) \Bigr) , \quad s \in [0,1]. 
$$

A point $(\xi,\tau) \in M_k \setminus B$ belongs to 
\begin{equation} \label{e:tildeen}
\left( \tilde e_n |_{M_k\setminus B} \right)^{-1} \bigl( F_n(B \times [0,1])\bigr)
\end{equation}
if and only if $\tilde e_n (\xi,\tau) = F_n(x,t,s)$ for some $(x,t,s) \in B \times [0,1]$.
Comparing the components in $\R^{k+1} \times \{0\}$ and $\{0\} \times \R^{k+1}$, this translates to
the two equations
\begin{eqnarray}
\xi + c \cos (\tau + 2\delta) \2 \xi &=& x+ c \cos (t+2\delta) \2 x , \\
\widetilde \sin \2 \tau \, w_0 &=& \widetilde \sin \2 t \, \bigl( s w_0 + (1-s) w_n(x) \bigr) . \label{e:2}
\end{eqnarray}
Let $(\xi,\tau) \in M_k \setminus B$ and $(x,t,s) \in B \times [0,1]$ be a solution.
Since $c<1$, the first equation shows that $\xi = x$ and $\cos (\tau+2\delta) = \cos (t+2\delta)$,
i.e., 
\begin{equation} \label{e:tt}
t+\tau + 4\delta \,=\, 2\pi .
\end{equation}
For $(x,t) \in B$ we have $x \in S_{\geq 0}^k := \left\{ x \in S^k \mid x_{k+1} \geq 0 \right\}$.
Therefore $\xi = x \in S_{\geq 0}^k$, whence $(\xi,\tau) \in (M_k\setminus B) \cap (S_{\geq 0}^k \times S^1)$,
see the shaded region in Figure~\ref{fig-}. 
Hence $\tau \in \, ]-\delta,\delta[$.
In particular, $\widetilde \sin \tau =0$.
By~\eqref{e:2},
$$
\widetilde \sin \2 t \, \bigl( s w_0 + (1-s) w_n(x) \bigr) \,=\, 0 .
$$
Assume that $\widetilde \sin \2 t = 0$. Then $t \in [-\gd,\gd]$ or $t=\pi$,
which are both impossible in view of~\eqref{e:tt} and the fact that $\delta$ is small.
Hence $\widetilde \sin \2 t \neq 0$, and so $s w_0 + (1-s) w_n(\xi) =0$,
i.e., 
$$
s=\tfrac 12 \quad \mbox{ and } \quad w_n(\xi) = - w_0 .
$$
By our choice of $w_0$ and $w_n$, the latter equation has exactly $n$ solutions
$\xi_1, \dots, \xi_n \in S^k$.
The closure of the set~\eqref{e:tildeen} is therefore
$$
\bigl\{ (\xi_j,t) \mid t \in [-\gd, \gd], \, j=1, \dots, n \bigr\} .
$$
The class represented by this set in $H_1 ( M_k \setminus \Int B, \pp B; \ZZ_{M_k}) \cong H_1(S^k \times S^1;\ZZ_{M_k}) \cong \ZZ_{M_k}$
is the one represented by the $n$~circles $\bigcup_{j=1}^n \{\xi_j\} \times S^1$.
The lemma is proved.
\proofend

\section{Where topology ends and geometry starts} 
There are several meaningful viewpoints on the phenomenon of Lagrangian (un)knots:
a surprising interplay between differential geometry and topology, 
an $h$-principle or its violation, 
a topological constraint on an invariant set of a dynamical system... 
We  conclude the paper with yet another aspect of Lagrangian {\it un}\1knottedness. 
The space~$\Lag$ of Lagrangian submanifolds which are Hamiltonian isotopic to a given one admits a transitive action of 
the group~$\Ham$ of Hamiltonian diffeomorphisms.  
Thus various interesting metric structures on~$\Ham$ that are known since the birth of symplectic topology -- such as Hofer's bi-invariant Finsler metric and Viterbo's spectral metric -- descend to~$\Lag$. 
The exploration of the corresponding geometry of~$\Lag$ became an active research area. 
A fair exposition would require a separate survey, so we invite the reader to google ``the Lagrangian Hofer metric". 
Here we just conclude with the following motto: 
{\it Lagrangian unknottedness paves a way for non-trivial measurements on the space of Lagrangian submanifolds.}

\m
\ni
{\bf Acknowledgment.}
We thank Yasha Eliashberg for exciting discussions, 
Jo\'e Brendel for the two proofs of Proposition~\ref{p:abc},
Kai Cieliebak for the proof of Proposition~\ref{p:GLh},
Andr\'as Juh\'asz for the proof of Proposition~\ref{p:Ju},
Stefan Nemirovski for the construction of exotic embeddings of $S^k \times S^1$ in~$\CC^{k+1}$,
and Denis Auroux, Jo\'e Brendel, Georgios Dimitroglou Rizell, Misha Entov, Jonny Evans, Joel Schmitz, and Arkadiy Skopenkov
for many useful remarks.
FS is deeply grateful to Yasha for arranging for him 
a wonderful sabbatical at Stanford in Spring~2024.


\section{Appendix by Georgios Dimitroglou Rizell}

The goal of this appendix is to show how well-established and classical results concerning the uniqueness of symplectic fillings of three-dimensional contact manifolds can be used to give short proofs 
of classification results for certain Lagrangian surfaces. More precisely, we show the following.

\begin{theorem} \label{thm:lagdisc}
Let $\Sigma$ be either $\R^2$ or $S^2$, and $L \subset (T^*\Sigma,d(p\,dq))$ 
a properly embedded open Lagrangian disc which coincides with a cotangent fibre $F=T^*_{\OP{pt}}\Sigma$ 
outside of a compact subset. 
Then there is a compactly supported symplectomorphism $\phi \in \OP{Symp}_c(T^*\Sigma)$ for which $\phi(F)=L$.
\end{theorem}

\begin{remarks*}
{\rm
(i) 
Gromov has shown that $\OP{Symp}_c(T^*\R^2)=\OP{Symp}_c(\R^4)$ is weakly contractible,~\cite{Gr85}. 
In particular, in this case, $L$ is compactly supported Hamiltonian isotopic to~$F$. 
This is a special case of the main result of Eliashberg and Polterovich from~\cite{EP96}; 
see Theorem~\ref{t:nolocal} in the article.

\s
(ii)
Seidel has shown that $\OP{Symp}_c(T^*S^2)\sim \{\tau^l\}_{l\in\Z}$ is weakly homotopy equivalent 
to the infinite cyclic group generated by the Dehn twist, \cite[Proposition 2.4]{Se08}. 
In particular, in this case, $L$ is compactly supported Hamiltonian isotopic to the image $\tau^l(F)$ of~$F$ 
under some power of the Dehn twist. 
To the author's knowledge, this result did not appear in the literature before. 
See the work \cite{CDR23} by C\^{o}t\'{e} and the author for similar results in the case when 
$\Sigma$ is an open Riemann surface.
}
\end{remarks*}

Before we prove the theorem we will give some background as well as a discussion about the method that we use. 
In~\cite{EGL20}, Eliashberg, Ganatra, and Lazarev introduced the concept of regular Lagrangians for studying the classification of certain Lagrangians up to symplectomorphism. One formulation of regularity is that the complement of the Lagrangian has the structure of a Weinstein cobordism. The question about which Lagrangians are regular is very hard, and little is known in general. In any case, assuming regularity, the classification problem of these Lagrangians can be reduced to understanding Weinstein handle decompositions of the complementary Weinstein cobordisms. This perspective was successfully used by Lazarev in~\cite{La20}.

Here we illustrate the power of the above perspective in the case of certain four-dimensional symplectic manifolds, where several strong uniqueness results are known for the symplectic structures on domains and cobordisms 
with fixed contact boundary. We will concentrate on the case when $L^n \subset (X^{2n},d\eta)$ is a properly embedded open Lagrangian disc inside a Liouville manifold. 
When talking about Liouville and Weinstein manifolds~$(X,d\eta)$ we will always assume that they are 
of finite type. 
We will also require that $L$ is cylindrical over a Legendrian outside of a compact subset, 
i.e.~that it is tangent to the Liouville vector field there; 
one calls $X$ a Liouville filling of its ideal contact boundary $(\partial_\infty X,\xi)$ at infinity, 
and $L$ a Lagrangian filling of its ideal Legendrian boundary $\partial_\infty L \subset \partial_\infty X$. 
Being of finite type is equivalent to the existence of a compact Liouville or 
Weinstein domain $\overline{X} \subset (X,d\eta)$ inside the manifold which has smooth boundary of contact type, 
i.e.~along which the Liouville flow is outwards transverse, and such that $(X \setminus \overline{X},d\eta)$ 
is a trivial Liouville cobordism, i.e.~symplectomorphic to
$$
\bigl( (0,+\infty)_t \times \partial \overline{X},d(e^t\alpha) \bigr), \quad \alpha=\eta|_{T\partial \overline{X}},
$$ 
while preserving the Liouville form. 
In particular, $(\partial\overline{X},\ker \alpha)$ is contactomorphic to $(\partial_\infty X,\xi)$.

It was shown in \cite[Proposition 2.3]{EGL20} that a Lagrangian disc inside a Liouville manifold is regular 
if and only if its complement can be given the structure of a Weinstein manifold; 
i.e.~the complementary Weinstein cobordism can be turned into a Weinstein manifold (and vice versa). 
This means that a Lagrangian disc is regular if and only if it arises as the Lagrangian co-core 
in a Weinstein handle decomposition of the Liouville manifold.

We elaborate a bit on the last statement. Recall the standard fact that, 
when $L \subset (X,d\eta)$ is an open and properly embedded Lagrangian disc in a Liouville manifold 
which is tangent to the Liouville vector field outside of a compact subset, 
then the complement $X \setminus L$ again admits the structure of a Liouville manifold $(X \setminus L,d \tilde \eta)$,
where $\tilde{\eta}=\eta+df$ is a suitable exact deformation by the differential of a function~$f$ 
supported near~$L$. 
The pair $L \subset (X,d\tilde{\eta})$ can then be recovered by a single Weinstein $n$-handle attachment 
along a Legendrian sphere $\Lambda' \subset \partial_\infty (X\setminus L)$, 
such that $L$ becomes the Lagrangian co-core of the corresponding handle attachment, 
and where $\Lambda=\partial_\infty L$ thus is the corresponding belt sphere.

Topologically the manifold $\partial_\infty(X \setminus L)$ is obtained by surgery on~$\partial_\infty X$ 
along the sphere~$\Lambda$. Further, as follows from the Lagrangian neighbourhood theorem applied to~$L$, the contact structure on $\partial_\infty(X \setminus L)$ is completely determined by 
the contact structure on~$\partial_\infty X$ together with 
the Legendrian sphere $\Lambda \subset \partial_\infty X$. 
This type of surgery can thus be performed in the category of contact manifolds along any 
Legendrian sphere $\Lambda \subset (Y,\xi)$, regardless of the existence of fillings 
of either the contact manifold or the Legendrian sphere; 
the resulting contact manifold is said to be obtained by $+1$-surgery along~$\Lambda$. 
This operation can be seen to be inverse to the usual contact $-1$-surgery that arises from standard 
Weinstein $n$-handle attachment; 
more precisely, performing a $+1$-surgery on the Legendrian belt-sphere 
that arises from the usual contact $-1$-surgery gives back the original contact manifold, and vice versa.

The main technical result used in the proof of Theorem \ref{thm:lagdisc} can be derived by techniques 
that go back to Gromov~\cite{Gr85} and Eliashberg~\cite{El90}; 
also Giroux' work~\cite{Gi91} can be used for classifying the contactomorphisms. 
For formulating the result, the following notion is crucial: 
A symplectomorphism 
$$
\bigl( \R_t \times Y_0,d(e^t\alpha_0) \bigr) \supset U \xrightarrow{\phi} 
\bigl( \R_t \times Y_1,d(e^t\alpha_1) \bigr)
$$
between two symplectisations will be called \textbf{cylindrical} if it is of the form 
$\phi(t,y)=(t-g(y),\psi(y))$. 
Consequently, we have $\psi^*\alpha_1=e^g\alpha_0$, i.e.~$\psi \colon (Y_0,\ker \alpha_0) \to (Y_1,\ker\alpha_1)$ 
is a contactomorphism. Moreover, the property of being cylindrical is equivalent to preserving the Liouville form.

\begin{proposition} \label{prp:uniqueness}
Let $(X,d\eta)$ be a Liouville manifold whose ideal contact boundary $\partial_\infty X$ is the standard tight contact~$S^3$, resp.\ $S^1 \times S^2$. 
Then
\begin{enumerate}
\item[{\rm (1)}]
The Liouville manifold $(X,d\eta)$ is symplectomorphic to 
$$
\left(\C^2,d\frac{1}{2}\sum_{i=1}^2 r_i^2\,d\theta_i=\omega_0\right),\ \:\:\text{resp.}\:\: \left(\C^*\times \C,d\left(\log{r_1}\,d\theta_1+ \frac{1}{2}r_2^2\,d\theta_2\right)\right),
$$
\item[{\rm (2)}]
Any cylindrical symplectomorphism
$$
\phi \colon X \setminus K \xrightarrow{\cong} X \setminus \phi(K),
$$
for some compact $K \subset X$, extends to a symplectomorphism of~$X$, 
possibly after first enlarging the compact subset~$K$.
\end{enumerate}
\end{proposition}

\begin{remark*}
{\rm
The uniqueness of Liouville fillings in (1) more generally holds for all contact three-manifolds that admit a subcritical Stein filling $X^2 \times \C$, 
see Theorem~16.9 in the book~\cite{CE12} by Cieliebak and Eliashberg. 
Further, it should be possible to use the same proof with more care to show that the uniqueness of filling 
holds relative to a fixed identification of the boundary, i.e.\ that (2) holds in these cases as well. 
Once this is established, our proof of Theorem~\ref{thm:lagdisc} immediately extends to fibers
$L \subset T^*\Sigma$ for an arbitrary surface~$\Sigma$ (open or closed).
}
\end{remark*}

\proof
\emph{The case when $\partial_\infty X=S^3$:}

(1): This is Gromov's result from \cite[0.3.C]{Gr85}.

(2): This folklore result was proven by Casals--Sp\'{a}\v{c}il in \cite[Theorem 4]{CS16}. 
Alternatively, one can verify that Gromov's construction of a symplectomorphism extends the given identification of the ideal contact boundary; i.e.~given an identification of the two fillings outside of a compact subset, the symplectomorphism that Gromov constructs can be taken to preserve this identification.

\s
\emph{The case when $\partial_\infty X=S^1 \times S^2$:}

(1): By \cite[Theorem 16.9]{CE12} the Liouville manifold~$X$ is obtained from a filling of 
the standard contact~$S^3$ by a standard Weinstein 1-handle attachment. 
The result then follows from the uniqueness of the Liouville filling of~$S^3$.

(2): This follows from Min's result \cite[Theorem 1.3]{Min24}, 
which builds heavily on work of Ding and Geiges~\cite{DG10}. 
It is shown that the contactomorphism group satisfies 
$\pi_0(\OP{Cont}(S^1 \times S^2,\xi_{std})) \cong \Z \times \Z_2$.

The generators $(1,0)$ and $(0,1) \in \Z \times \Z_2$ can be seen to  be induced by the following symplectomorphisms. First, $(0,1)$ is induced by the holomorphic involution
$$
\begin{array}{rcl}
\C^* \times \C &\to& \C^* \times \C,\\
(z_1,z_2) &\mapsto& (z_1^{-1},z_2),
\end{array}
$$
which preserves the Liouville form. The generator $(1,0)$ is constructed from a Dehn twist on a two-torus 
inside $S^1 \times S^2$, see \cite[Section 3]{DG10}. 
This contactomorphism is induced by 
$$
\begin{array}{rcl}
\C^* \times \C &\to& \C^* \times \C,\\
(r_1,r_2,\theta_1,\theta_2) &\mapsto& \left(e^{-\frac{1}{2}r_2^2}r_1,r_2,\theta_1,\theta_2+\theta_1\right),
\end{array}
$$
(i.e.~a symplectic suspension of a full $2\pi$-rotation of $\C$) where we have used polar coordinates 
on each factor for the description. 
Note that this symplectomorphism also preserves the Liouville form.
\proofend

\ni
{\it Proof of Theorem \ref{thm:lagdisc}.}
Consider the Lagrangian disc $L \subset T^*\Sigma$. 
In the case when $\Sigma=\R^2$ we can find an exact deformation of the tautological Liouville form 
after which $T^*\R^2$ becomes identified with the Liouville manifold 
$\left(\C^2,d\frac{1}{2}\sum_i(x_i\,dy_i-y_i\,dx_i)\right)$.

We thus consider the two cases when either $X=\C^2$ or $X=T^*S^2$, 
and denote by $\Lambda := \partial_\infty L \subset \partial_\infty X$ the ideal Legendrian boundary 
of the Lagrangian disc $L \subset X$. 
In the case $X=\C^2$ it follows that $\Lambda$ is the standard Legendrian unknot inside the standard 
contact sphere $\partial_\infty X=S^3$, 
while in the case $X=T^*S^2$ it is the Legendrian spherical cotangent fibre.

In both cases the ideal contact boundary produced by a contact $+1$-surgery performed along 
$\Lambda \subset \partial_\infty X$ is easy to find explicitly. 
Indeed, if we remove a cotangent fibre from~$T^*\R^2$, we obtain $T^*(S^1 \times \R )=\C^* \times \C$ 
with ideal contact boundary given by the standard tight $\partial_\infty (\C^* \times \C)=S^1 \times S^2$; 
if we remove a cotangent fibre from~$X=T^*S^2$, we obtain $T^*\R^2=\C^2$ with ideal contact boundary 
$\partial_\infty \C^2=S^3$, i.e.~the standard tight contact three-sphere.

Consider a sufficiently large compact subset $\overline{X} \subset X$. 
Using the symplectic standard neighbourhood theorem, we can readily construct a symplectomorphism
$$
\phi \colon (X \setminus \overline{X}) \cup \mathcal{O}(F) \xrightarrow{\cong} (X \setminus \overline{X}) \cup \mathcal{O}(L) 
$$
where $\mathcal{O}(F)$ and $\mathcal{O}(L)$ denote open neighbourhoods of $F$ and~$L$, respectively, 
and such that $\phi(F)=L$. Taking additional care, we may further assume that
$$
\phi|_{X \setminus \overline{X}}=\id_{X \setminus \overline{X}}
$$
is satisfied. This enables us to choose a Liouville form on $X \setminus L$ which, 
outside of a compact subset, is defined uniquely by the requirement that it pulls back to the Liouville form 
on $X \setminus F$ under 
$\phi|_{\left( (X \setminus \overline{X}) \cup \mathcal{O}(F) \right) \setminus F}$.

In particular, the symplectomorphism $\phi|_{(X \setminus F)\setminus K}$ is cylindrical for some sufficiently large compact $K \subset X \setminus F$. 
Proposition~\ref{prp:uniqueness} then provides the sought extension of $\phi$ to a compactly supported symplectomorphism $(X,F) \xrightarrow{\cong} (X,L)$.
\proofend

\end{document}